\documentclass[11pt]{article}
\usepackage{amssymb}
\usepackage{amsmath}
\usepackage{amscd}
\usepackage{amsfonts}
\usepackage{graphicx,color}
\usepackage[title]{appendix}
\usepackage{mathtools}

\numberwithin{equation}{section}

\newtheorem{theorem}{Theorem}[section]
\newtheorem{corollary}{Corollary}[section]
\newtheorem{lemma}{Lemma}[section]
\newtheorem{remark}{Remark}[section]

\makeatletter
\newcommand{\vast}{\bBigg@{3}}
\newcommand{\Vast}{\bBigg@{4}}

\begin{document}
\title{\bf On the proof of a variant of Lindel{\"o}f's hypothesis}
\author{by \\\\  A.S. Fokas}

\date{}
\maketitle
\begin{center}
Department of Applied Mathematics and Theoretical Physics,

University of Cambridge, CB3 0WA, UK,

and

Viterbi School of Engineering, University of Southern California, 

Los Angeles, California, 90089-2560, USA. 
\end{center}






\begin{abstract}
The leading asymptotic behaviour as $t\to \infty$ of the celebrated Riemann zeta function $\zeta(s), \ s = \sigma + it,  \quad 0<\sigma<1, \quad t>0 , \ t\to\infty,$ can be expressed in terms of a transcendental sum. The sharp estimation of this sum remains one of the most important open problems in mathematics with a long and illustrious history. Lindel{\"o}f's hypothesis states that for $\sigma=1/2$, this sum is of order $O\left(t^\varepsilon\right)$ for every $\varepsilon>0$. We have recently introduced a novel approach for estimating such transcendental sums: we first embed the Riemann zeta function in a more complicated mathematical structure, and then compute the large $t$-asymptotics of this structure. In particular, we have embedded the Riemann zeta function in a certain Riemann-Hilbert problem and we have began the analysis of the large $t$-asymptotics of the associated integral equation. We have shown that the occurrence of certain gamma functions in the above integral equation motivates splitting the relevant interval of integration into three subintervals which are defined in terms of the small positive numbers $\delta_1$ and $\delta_4$. The asymptotic analysis of the resulting integral equation requires the further splitting of the relevant interval of integration into four subintervals which are defined in terms of the small positive numbers $\{\delta_j\}_1^4$. The rigorous asymptotic analysis of the first two relevant integrals, $I_1$ and $I_2$, was performed in \cite{F}. Here, the rigorous analysis is performed of the last two integrals, $I_3$ and $I_4$. The combination of the above results yields a proof for the analogue of  Lindel{\"o}f's hypothesis for a slight variant of the transcendental sum characterising the large $t$-asymptotics of $|\zeta(s)|^2$, namely for a sum which differs from the latter sum only in the occurrence of a logarithmic term which is larger than $\frac{1}{2}\ln t$ and smaller than $t^{\varepsilon}$. Interestingly, the parameter $\varepsilon$ in  Lindel{\"o}f's hypothesis is explicitly defined in terms of $\delta_3$. 
\end{abstract}

\section{Introduction}

The Riemann zeta function occurs in many areas of mathematics, and in particular it plays central role in analytic number theory. Several conjectures related to the Riemann function remain open, including the Riemann hypothesis, perhaps the most celebrated open problem in the history of mathematics. The latter hypothesis states that

\begin{equation*}
\zeta(s) \ne 0, \quad s = \sigma + it, \quad \text{for} \quad 0<\sigma<\frac{1}{2}, \quad t \in \mathbb{R}.
\end{equation*}
This hypothesis can be verified numerically for $t$ up to order $O(10^{13})$, thus the basic problem associated with the Riemann hypothesis is its proof for large $t$. This provides an additional motivation for studying the asymptotic behavior of $\zeta(s)$ as $t \to \infty$. This problem, which has a long and illustrious history, is deeply related with the Lindel{\"o}f hypothesis. Indeed, it is well known that the leading asymptotic behavior as $t \to \infty$ of $\zeta(s)$, $0<\sigma<1$, $t>0$, can be expressed in terms of the following transcendental sum:

\begin{equation} \label{1.1}
\zeta(s) \sim \sum_{m=1}^{[T]} \frac{1}{m^s}, \quad T = \frac{t}{2\pi}, \quad t \to \infty,
\end{equation}
where throughout this paper $[A]$ denotes the integer part of the positive number $A$.

The estimation of this sum remains a most important open problem. Lindel{\"o}f's conjecture states that for $\sigma=\frac{1}{2}$ this sum is of order $O(t^{\varepsilon})$ for any $\varepsilon>0$. The Riemann hypothesis implies Lindel{\"o}f hypothesis, and conversely, Lindel{\"o}f's hypothesis implies that very ``few zeros can escape Riemann's hypothesis" \cite{GM}. Indeed, let $N(\sigma,t)$ denote the number of zeros, $\beta+it$, of $\zeta(s)$, such that $\beta>\sigma$ and $0<t\leq T$. It is stated in \cite{Tu} that ``Lindel{\"o}f's hypothesis is much stronger than expected and even implies the estimate
\begin{equation*}
N(\sigma,t) = O(T^{\varepsilon}), \quad \frac{1}{2} + \delta \leq \sigma <1, \quad T \to\infty,
\end{equation*}
for $\varepsilon$ and $\delta$ positive and arbitrarily small''.

The sum of the rhs of \eqref{1.1} is a particular case of an exponential sum. Pioneering results for the estimation of such sums were obtained almost 100 years ago using methods developed by Weyl \cite{W}, and Hardy and Littlewood \cite{HL}, when it was shown that $\zeta(1/2 + it) = O(t^{1/6 +  \varepsilon})$. In the last 90 years some slight progress was made using the ingenious techniques of Vinogradov \cite{V}. 
Further progress was made by several authors, and currently the best results is due to Bourgain \cite{B},  who was able to reduce the exponent factor to $53/342\approx 0.155$.

Regarding the large $t$ asymptotics of Riemann's zeta function, we note that the best estimate for the growth of $\zeta(s)$ as $t \to \infty$ is based on the approximate functional equation, see page 79 of \cite{T},

\begin{multline} \label{1.2}
\zeta(s) = \sum_{n\le x} \frac{1}{n^s} + \frac{(2\pi)^s}{\pi} \sin{\left(  \frac{\pi s}{2} \right)} \Gamma (1-s) \sum_{n\le y} \frac{1}{n^{1-s}}  + O \left(  x^{-\sigma} + |t|^{\frac{1}{2}-\sigma} y^{\sigma -1}  \right),    \\
xy = \frac{t}{2\pi}, \quad 0<\sigma<1, \quad t \to \infty,
\end{multline}
where $\Gamma(s)$, $s \in \mathbb{C}$, denotes the gamma function. It should be emphasized that, in contrast to the usual situation in asymptotics, where higher order terms in an asymptotic expansion are more complicated, the higher order terms of the asymptotic expansion of  $\zeta(s)$ can be computed {\it explicitly}. Siegel, in his classical paper \cite{S}, building on Riemann's unpublished notes, presented the asymptotic expansion of $\zeta(s)$ to {\it all} orders in the important case of $x=y=\sqrt{t/2\pi}$. In \cite{FL}, analogous results are presented for {\it any} $x$ and $y$ valid to {\it all} orders; these results play a crucial role in the asymptotic analysis presented here (similar results for the Hurwitz function are presented in \cite{FF}).

We have recently introduced a new approach for estimating the large $t$ asymptotics of the Riemann zeta and related functions: instead of analyzing $\zeta(s)$ directly, we first obtain a particular equation satisfied by $\zeta(s)$, and then compute the large $t$ asymptotics of this equation. Actually, it is shown in \cite{F} that $\zeta(s)$ satisfies a certain Riemann-Hilbert problem, or equivalently, the following singular integral equation:

\begin{multline} \label{1.3}
\frac{t}{\pi}  \oint_{-\infty}^{\infty} \Re{\left\{   \frac{\Gamma(it - i\tau t)}{\Gamma(\sigma + i t)} \Gamma(\sigma + i\tau t) \right\}} \left| \zeta (\sigma+i\tau t) \right|^2 \text{d}\tau + \mathcal{G}(\sigma,t) = 0, \\
0<\sigma< 1, \quad t>0,
\end{multline}
where the principal value integral is defined with respect to $\tau=1$ and the function $\mathcal{G}(\sigma,t)$ is defined by

\begin{multline} \label{1.4}
\mathcal{G}(\sigma,t) = \begin{cases}
      \zeta(2\sigma) + \left( \frac{\Gamma(1-\bar{s})}{\Gamma(s)} + \frac{\Gamma(1-s)}{\Gamma(\bar{s})} \right) \Gamma(2\sigma-1) \zeta(2\sigma-1) + \frac{2(\sigma -1)\zeta(2\sigma-1)}{(\sigma - 1)^2 + t^2}, & \sigma \ne \frac{1}{2}, \\
     \Re{\left\{ \Psi \left(  \frac{1}{2} + it \right) \right\}} + 2\gamma - \ln{2\pi} + \frac{2}{1+4t^2}, & \sigma = \frac{1}{2},      
   \end{cases}
\end{multline}
with $\Psi(z)$ denoting the digamma function, i.e.,

\begin{equation*}
\Psi(z) = \frac{\frac{\mathrm{d}}{\mathrm{dz}}\Gamma(z)}{\Gamma(z)}, \quad z \in \mathbb{C},
\end{equation*}
and $\gamma$ denoting the Euler constant.

We recall that the theory and applications of Riemann-Hilbert problems \cite{I} have flourished in the last 40 years, with spectacular applications from the theory of integrable systems \cite{Z} to the asymptotics of orthogonal polynomials \cite{FIK}.

The term $ \Gamma(it - i\tau t)  \Gamma(\sigma + i t) / \Gamma(\sigma + i\tau t)$, $-\infty<\tau<\infty$, occurring in \eqref{1.3}, decays exponentially for large $t$, unless
 \begin{equation}\label{1.5new}
 -t^{\delta_{1}-1}\leq \tau\leq 1+t^{\delta_{4}-1},
 \end{equation}
 where $\delta_{1}$ and $\delta_{4}$ are sufficiently small positive constants. Thus, equation \eqref{1.3} simplifies to the equation (for rigorous details see \cite{F})
\begin{equation} \label{1.5}
\frac{t}{\pi} \oint_{-t^{\delta_1-1}}^{1+t^{\delta_4-1}} \Re {\left\{ \frac{\Gamma(it-i\tau t)}{\Gamma(\sigma + it)}\Gamma(\sigma + i\tau t) \right\}} \left| \zeta(\sigma+i\tau t)\right|^2 \textrm{d}\tau + \mathcal{G}(\sigma,t)  \nonumber
\end{equation}
\begin{equation}
+ O\left(e^{-\pi t^{\tilde{\delta}_{14}}}\right) = 0, \quad 0<\sigma<1, \quad t\to\infty,
\end{equation}
where the principal value integral is defined with respect to
$\tau=1$, and $\tilde{\delta}_{14}=\mbox{min}(\delta_{1},\delta_{4})$.

The computation of the large $t$ asymptotics of \eqref{1.5}
requires splitting further the interval $[-t^{\delta_{1}-1},1+t^{\delta_{4}-1}]$
into the following four subintervals:
\begin{align}\label{1.6}
&L_{1}=[-t^{\delta_{1}-1},t^{-1}],
L_{2}=[t^{-1},t^{\delta_{2}-1}], L_{3}=[t^{\delta_{2}-1},1-t^{\delta_{3}-1}],\notag\\
&L_{4}=[1-t^{\delta_{3}-1},1+t^{\delta_{4}-1}],
\end{align}
where $\delta_{2}$ and $\delta_{3}$ are sufficiently small positive  constants. Thus, the asymptotic evaluation of \eqref{1.5} reduces to the analysis of the four integrals,
\begin{equation}
I_{j}(\sigma,t)=\frac{t}{\pi} \oint_{L_{j}} \Re {\left\{ \frac{\Gamma(it-i\tau t)}{\Gamma(\sigma + it)}\Gamma(\sigma + i\tau t) \right\}} \left| \zeta(\sigma+i\tau t)\right|^2 \textrm{d}\tau ,
\end{equation}
where $I_{1}$, $I_{2}$, $I_{3}$, $I_{4}$ also depend on $\delta_{1}$,
$\delta_{2}$, $(\delta_{2},\delta_{3})$, $(\delta_{3},\delta_{4})$, respectively,  $0<\sigma<1$, $t>0$,  $\{L_{j}\}_{1}^{4}$
are defined in \eqref{1.6}, and the principal value integral is needed only for $I_{4}$.

The rigorous estimation of $I_{1}$ and $I_{2}$ was performed in \cite{F}, and it was shown that these integrals are ``small". Using the precise estimates of $I_1$ and $I_2$ presened in \cite{F}, as well as the large $t$-asymptotics of $\mathcal {G}(\sigma,t)$, equation \eqref{1.5} yields
\addtocounter{equation}{1}
\begin{multline} \label{4.21}
I_{3}(\sigma ,t,\delta _{2},\delta _{3}) + I_4(\sigma, t,\delta_3,\delta_4)  +     \ln{t} + 2\gamma - \ln{2\pi} + O \left(  t^{\delta_2-\frac{1}{2}} \ln{t} \right),\\  
   + O \left(  e^{-\pi t^{\tilde{\delta}_{14}}} \right) + O\left(  \frac{t^{\delta_1}}{t^{\sigma}} \right) + O \left( \frac{1}{t^2}\right) = 0, \qquad \sigma = \frac{1}{2}, \quad t\to\infty, \tag{1.9a}
\end{multline}
and
\begin{multline} \label{4.22}
I_{3}(\sigma ,t,\delta _{2},\delta _{3}) + I_4(\sigma, t,\delta_3,\delta_4)+ \zeta(2\sigma) + 2\Gamma(2\sigma - 1) \zeta(2\sigma - 1) \sin{(\pi\sigma) t^{1-2\sigma}} \left(  1+O\left( \frac{1}{t}  \right) \right) \\
+ \begin{cases}
     O\left( t^{-\sigma + \left(\frac{3}{2} -\sigma  \right)\delta_2} \ln{t} \right), & 0<\sigma < \frac{1}{2} \\
     O\left( t^{-\sigma + \left(\sigma + \frac{1}{2} \right)\delta_2}\zeta(2\sigma)\right), &  \frac{1}{2} < \sigma <1   
   \end{cases}\\
   + O \left(  e^{-\pi t^{\tilde{\delta}_{14}}} \right) + O\left(  \frac{t^{\delta_1}}{t^{\sigma}} \right) + O \left( \frac{1}{t^2}\right) = 0, \quad t\to\infty. \tag{1.9b}
\end{multline}

It was shown in \cite{F} that if one replaces the term $|\zeta(\sigma+it)|^2$ occurring in the integral $I_3$ by the leading asymptotics of $|\zeta(\sigma+i\tau t)|^2$ as $\tau t\to\infty$, then the leading asymptotics of $I_3$ involves two contributions: the contribution from the associated stationary points, denoted by $I_{S}$, and the contribution from the neighbourhood of the point $1-t^{\delta_{3}-1}$, which is denoted by $I_{B}$. The former contribution, to the leading order, is given by
\begin{equation} \label{1.7}
I_S \sim 2 \Re{\left\{  \mathop{\sum\sum}_{m_1,m_2 \in M}  \frac{1}{m_1^s (m_1+m_2)^{\bar{s}}}  \right\}}, \quad 0<\sigma<1, \quad t\to\infty,
\end{equation}
where the set $M$ is defined by
\begin{multline}  \label{1.8}
M = \left\{ m_{1}\in\mathbb{N}^{+}, \quad  m_{2}\in\mathbb{N}^{+}, \quad 1\leq m_{1}\leq [T], \quad 1\leq m_{2}< [T], \right.  \\
\left. \frac{1}{t^{1-\delta _{3}}-1}<\frac{m_{2}}{m_{1}}<t^{1-\delta _{2}}-1, \quad  t>0, \quad  T=\frac{t}{2\pi} \right\}.
\end{multline}

The latter contribution to the leading order is given by
\begin{subequations} \label{112all}
\begin{align} \label{1.9new}
I_B (\sigma, t, \delta_3) \sim & t^{-\frac{\delta_3}{2}}  t^{i(\delta_{3}-1)t^{\delta_3}}(1-t^{\delta_{3}-1})^{\sigma-\frac{1}{2}+
i(t-t^{\delta_{3}})}\notag\\
&\times \mathop{\sum\sum}_{m_1,m_2 \in N}\frac{1}{m_1^{s-it^{\delta_{3}}}}
\frac{1}{m_2^{\bar{s}+it^{\delta_{3}}}}\frac{1}{
\ln\left[\frac{m_2}{m_1}(t^{1-\delta_3}-1)\right]} , 
\end{align}
where the set $N$ is defined by
\begin{multline}  \label{1dN}
N= \bigg\{ m_{1}\in\mathbb{N}^{+}, \quad  m_{2}\in\mathbb{N}^{+}, \quad 1\leq m_{1}\leq [T], \quad 1\leq m_{2}< [T],  \\
 \frac{m_{2}}{m_{1}}>\frac{1}{t^{1-\delta _{3}}-1}\bigg(1+c(t)\bigg), \ t^{-\frac{\delta_3}{2}}\ll c(t) \ll 1, \quad  t>0, \quad  T=\frac{t}{2\pi} \bigg\}. 
\end{multline}
\end{subequations}

Similarly, it was shown in \cite{F} that if one replaces the term $|\zeta(\sigma+it)|^2$ occurring in the integral $I_4$ by the leading order asymptotics of $|\zeta(\sigma+i\tau t)|^2$ as $\tau t\to\infty$, then the leading asymptotics of $(-I_{4})$
is the same with the leading asymptotics of $|\zeta(s)|^{2}$, namely,
\begin{equation} \label{1.10}
I_4 \sim - \sum_{m_1=1}^{[T]}  \sum_{m_2=1}^{[T]} \frac{1}{m_1^s m_2^{\bar{s}}}, \quad 0<\sigma<1, \quad t\to\infty.
\end{equation}
%
%
%

It is remarkable that  for large $t$, the sums of the rhs of equations \eqref{1.7} and \eqref{1.10} are identical. The starting point for proving this result is the derivation in \cite{F} of two exact identities, given below 
\begin{equation}\label{1.13-new}
2\Re\left\{
\sum_{m_{1}=1}^{[T]}\sum_{m_{2}=1}^{[T]}\frac{1}{m_{1}^{s}
(m_{1}+m_{2})^{\bar{s} }} \right\} =
 \sum_{m_{1}=1}^{[T]} \sum_{m_{2}=1}^{[T]}\frac{1}{m_{1}^{s}m_{2}^{\bar{s} }} \nonumber
\end{equation}
\begin{equation}  \label{sevensix}
-\sum_{m=1}^{[T]}\frac{1}{m^{2\sigma }}
+2\Re \left\{ \sum_{m=1}^{[T]}\sum_{n=[T]+1}^{[T]+m}\frac{1}{m^{\bar{s} }n^{s}} \right\},
\hphantom{3a} s\in\mathbb{C},
\end{equation}
and
\begin{equation} \label{1.11}
\sum_{m_1=1}^{[T]}  \sum_{m_2=1}^{[T]} \frac{1}{m_1^s (m_{1}+m_2)^{\bar{s}}}=
 \sum_{m_1,m_2\in M}\frac{1}{m_1^s (m_{1}+m_2)^{\bar{s}}}+S_{1}+S_{2},
\end{equation}
where the sum $S_{1}$ and $S_{2}$ are defined by
\addtocounter{equation}{1}
\begin{equation}\label{1.12}
S_{1}(\sigma,t,\delta )=\sum_{m_{1}=1}^{\frac{[T]}{t^{1-\delta }-1}-1}\sum_{m_{2}=(t^{1-\delta }-1)m_{1}+1}^{[T]}
\frac{1}{m_{1}^{s}(m_{1}+m_{2})^{\bar{s} }},\tag{\theequation a}
\end{equation}
and
\begin{align} \label{1.13}
S_{2}(\sigma,t,\delta )=\sum_{m_{1}=t^{1-\delta}}^{[T]} \sum_{m_{2}=1}^{\frac{m_1}{t^{1-\delta}-1}-1}
\frac{1}{m_{1}^{s}(m_{1}+m_{2})^{\bar{s} }},\tag{\theequation b}
\end{align}
with $s=\sigma+it, ~ 0<\sigma<1, ~t>0$.

In this paper the following results are derived: first, the next term of the large $\tau t$ asymptotics of $|\zeta(\sigma+i\tau t)|^2$ is included in $I_3$ and $I_4$, and it is shown that this term contributes terms in $I_{3}+I_{4}$ which are of order $O(t^{\frac{\delta}{2}}\ln t)$. Second, the sums $S_{1}$ plus $S_{2}$ plus the last two
terms of the rhs of \eqref{1.13-new} are of order $O(t^{\frac{\delta}{2}}\ln t)$
for $\sigma=\frac{1}{2}$, and   of  order $O(1)$ for $\frac{1}{2}<\sigma<1$.

Thus, combining equations (1.8)-(1.15) with the rigorous estimates presented here, we conclude that
\begin{multline} \label{1.14}
\mathop{\sum\sum}_{m_1,m_2 \in N} \frac{1}{m_1^{s-it^{\delta_{3}}}} \frac{1}{m_2^{\bar{s}+it^{\delta_{3}}}} \frac{1}{\ln{\left[ \frac{m_2}{m_1} \left( t^{1-\delta_{3}} - 1 \right)  \right]}} = \begin{cases} O \left(  t^{ \delta_{3} } \ln{t} \right), \ & \sigma=\frac{1}{2}, \\
O \left(  t^{ \frac{\delta_{3}}{2} }  \right), \  & \frac{1}{2}<\sigma<1,\end{cases}
   \\
 \quad t \to \infty,
\end{multline}
where $N$ is defined in \eqref{1dN}.

The main difference between the function defined by the lhs of \eqref{1.14} and the leading term of the large $t$-asymptotics of $|\zeta(s)|^{2}$ given by the sum of the
rhs of \eqref{1.10}, is the occurrence of the $\ln$ term, which satisfies
\begin{equation*}
\frac{1}{\ln t}<\frac{1}{\ln \left[\frac{m_2}{m_1}(t^{1-\delta_{3}}-1)\right]}<t^{\frac{\delta_{3}}{2}}.
\end{equation*}
Thus, equation \eqref{1.14} provides the analogue of Lindel\"{o}f's hypothesis
for a slight variant of $|\zeta(s)|^{2}$.

It is interesting to note that the occurrence of $\varepsilon$ in Lindel{\"o}f's hypothesis is somewhat mysterious, since $\zeta(1/2 + it)$ does {\it not} contain $\varepsilon$. On the other hand, the occurrence of the analogous small parameter $\delta_3$ appearing in the rhs of \eqref{1.14} is clearly explained: the asymptotic evaluation of \eqref{1.3} requires the introduction of the four small parameters $\{  \delta_j \}_{1}^{4}$, and at the end of a lengthy analysis, one of these parameters appears in the asymptotics of the relevant double sum.

A numerical comparison performed in \cite{F} of $|\zeta(s)|^{2}$ with the function defined
by the lhs of equation \eqref{1.14} suggest that for large $t$ the latter function times $\ln t$
can be approximated by $|\zeta(s)|^{2}$.
A possible approach for establishing rigorously such a relation between
the above two functions, and hence proving  Lindel\"{o}f's hypothesis, is discussed in \cite{FK}.

In the next two sections we present the asymptotic evaluation for large $t$ of $I_3$ and $I_4$. The analysis of the latter integral is easier than the former integral, thus in the next section we consider $I_4$.

\section{The Evaluation of $I_4$} \label{sec5}

\begin{theorem}\label{t5.1}
Let $I_4(\sigma, t,\delta_3,\delta_4)$ be defined by
\begin{equation} \label{5.1}
I_4(\sigma, t,\delta_3,\delta_4) = \frac{t}{\pi} \oint_{1-t^{\delta_3-1}}^{1+t^{\delta_4-1}}  \Re{\left\{ \Gamma(it - i\tau t)  \frac{\Gamma(\sigma + i\tau t)}{\Gamma(\sigma + i t)}\right\}} \left| \zeta(\sigma+i\tau t)\right|^2 \textrm{d}\tau ,   \nonumber
\end{equation}
\begin{equation}
 \quad 0<\sigma<1, \quad t>0,
\end{equation}
where $\zeta(z)$ and $\Gamma(z)$ denote the Riemann zeta and gamma functions respectively, $\delta_3$ and $\delta_4$ are sufficiently small, positive constants, and the principal value integral is defined with respect to $\tau=1$, i.e.,
\begin{equation} \label{5.2}
\oint_{1-t^{\delta_3-1}}^{1+t^{\delta_4-1}}  \textrm{d}\tau = \lim_{\varepsilon\to 0} \left( \int_{1-t^{\delta_3-1}}^{1-\varepsilon}  \textrm{d}\tau + \int_{1+\varepsilon}^{1+t^{\delta_4-1}}  \textrm{d}\tau  \right).
\end{equation}
Then,
\begin{equation} \label{5.3n}
I_4= - \sum_{m_1=1}^{[T]}  \sum_{m_2=1}^{[T]}\frac{1}{m_1^s m_2^{\bar{s}}}  + O\left(t^{1-2\sigma+\delta_{34}}\right),
\end{equation}
where
\begin{equation*}
T=\frac{t}{2\pi}, \quad \delta_{34} = \max{(\delta_3,\delta_4)}.
\end{equation*}
\end{theorem}

\textbf{Proof}
Using the change of variables $1-\tau=\rho$, $I_4$ becomes
\begin{equation*}
I_4 = \frac{t}{\pi} \oint_{-t^{\delta_4-1}}^{t^{\delta_3-1}} \Re{\left\{ \Gamma(it\rho)  \frac{\Gamma(\sigma + i t -it\rho)}{\Gamma(\sigma + i t)}\right\}} \left| \zeta(\sigma+i t-it\rho)\right|^2 \textrm{d}\rho,
\end{equation*}
where now the principal value integral is defined with respect to $\rho=0$. The change of variables $t\rho=x$ yields
\begin{equation} \label{5.4}
I_4 = \frac{1}{\pi} \oint_{-t^{\delta_4}}^{t^{\delta_3}} \Re{\left\{ \Gamma(ix)  \frac{\Gamma(\sigma + i t -ix)}{\Gamma(\sigma + i t)}\right\}} \left| \zeta(\sigma+i t-ix)\right|^2 \textrm{d}x.
\end{equation}
Since $|x|<t^{\delta_{34}}$, it follows that
\begin{equation} \label{5.5}
t-x\to\infty \quad \textrm{as} \quad t\to\infty,
\end{equation}
thus, we can use the large $\xi$ asymptotics of $\Gamma(\sigma+i\xi)$ with either $\xi=t-x$ or $\xi=t$ to compute the ratio $\Gamma(\sigma +it -ix) / \Gamma(\sigma + i t)$. In this connection we note that starting with Sterling's classical result, the following formulae are derived in \cite{FL}:

\begin{subequations} \label{219}
\begin{equation} \label{2.19a}
\Gamma(\sigma+i\xi) = \sqrt{2\pi}\xi^{\sigma-\frac{1}{2}}e^{-\frac{\pi \xi}{2}}e^{-\frac{i\pi}{4}}e^{-i\xi}\xi^{i\xi}e^{\frac{i\pi \sigma}{2}} \left[ 1 + O\left(\frac{1}{\xi}\right) \right], \quad \xi \to\infty, 
\end{equation}
and
\begin{equation} \label{2.19b}
\Gamma(\sigma-i\xi) = \sqrt{2\pi}\xi^{\sigma-\frac{1}{2}}e^{-\frac{\pi \xi}{2}}e^{\frac{i\pi}{4}}e^{i\xi}\xi^{-i\xi}
e^{-\frac{i\pi \sigma}{2}} \left[ 1 + O\left(\frac{1}{\xi}\right) \right], \quad \xi \to\infty. 
\end{equation}
\end{subequations}

\begin{figure}
\begin{center}
\includegraphics[width=.5\textwidth]{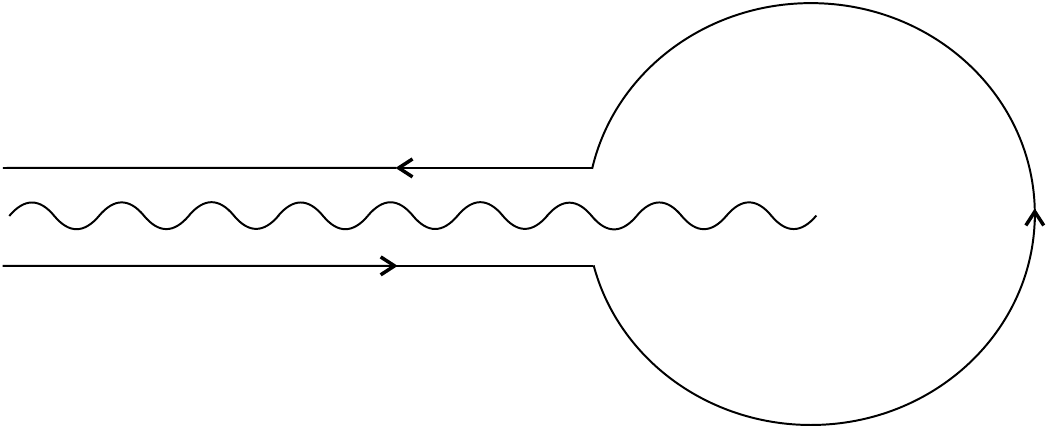}
\end{center}
\caption{The Hankel contour $H$.}
\label{fig1}
\end{figure}

Using the above formulae we find
\begin{equation}
\frac{\Gamma(\sigma +it -ix)}{\Gamma(\sigma + i t)} =
\nonumber
\end{equation}
\begin{equation}
\frac{t^{i(t-x)}\left(1-\frac{x}{t}\right)^{i(t-x)} t^{\sigma - \frac{1}{2}} \left(1-\frac{x}{t}\right)^{\sigma - \frac{1}{2}} e^{-\frac{\pi}{2}(t-x)} e^{-i(t-x)} e^{\frac{i\pi\sigma}{2}} \left[ 1 + O\left( \frac{1}{t-x} \right) \right]} {t^{it}t^{\sigma-\frac{1}{2}}e^{-it}e^{\frac{i\pi\sigma}{2}}e^{-\frac{\pi t}{2}}\left[ 1 + O\left( \frac{1}{t} \right) \right]} 
\nonumber
\end{equation}
\begin{equation}
=t^{-ix}  e^{\frac{\pi x}{2}} \left(1-\frac{x}{t}\right)^{\sigma - \frac{1}{2}} e^{ix} \left( 1 - \frac{x}{t} \right)^{i(t-x)} \left[ 1 + O\left( \frac{1}{t-x} \right) \right]
\nonumber
\end{equation}
\begin{equation}
\times \left[ 1 + O\left( \frac{1}{t} \right) \right] , \quad t\to\infty.  \nonumber
\end{equation}
Substituting the above expression in \eqref{5.4} we find
\begin{equation} \label{5.6}
I_4 = \frac{1}{\pi} PV\int_{-t^{\delta_4}}^{t^{\delta_3}} \Re{\left\{ \Gamma(ix) t^{-ix} e^{\frac{\pi x}{2}} e^{ix} \left(  1 - \frac{x}{t} \right)^{i(t-x)} \right\}} \left| \zeta(\sigma+i t-ix)\right|^2 
\nonumber
\end{equation}
\begin{equation}
\times \left[ 1 + O\left( \frac{1}{t-x} \right) \right] \left( 1-\frac{x}{t} \right)^{\sigma-\frac{1}{2}} \textrm{d}x \quad  
\left[ 1 + O\left( \frac{1}{t} \right) \right], \quad t\to \infty.
\end{equation}
In order to evaluate the above expression we will employ the representation
\begin{equation}
\Gamma(s) = \frac{1}{2i\sin{(\pi s)}} \int_{H} e^{z}z^{s-1}\textrm{d}z, \quad s\ne 0,-1,-2, \dots,
\end{equation}
where $H$ denotes the Hankel contour with a branch cut along the negative real axis, see figure \ref{fig1}, defined by
\begin{equation*}
H=\left\{ re^{-i\pi} ~ \Big| ~ 1<r<\infty \right\} \cup \left\{ e^{i\theta} ~ \Big| ~ -\pi<\theta<\pi \right\} \cup \left\{ re^{i\pi} ~ \Big| ~ 1<r<\infty \right\}.
\end{equation*}

\noindent Hence,
\begin{equation} \label{5.9}
\Gamma(ix) = \frac{1}{e^{-\pi x}-e^{\pi x}} \int_{H} \frac{e^z}{z}z^{ix}\textrm{d}z.
\end{equation}

It is shown in equation (3.2) of \cite{FL} that
\begin{equation} \label{5.10}
\zeta(\sigma + i\xi) = \sum_{m=1}^{\left \lfloor{\frac{\eta}{2\pi}}\right \rfloor}  \frac{1}{m^{\sigma+i\xi}} - \frac{1}{1-\sigma - i\xi} \left( \frac{\eta}{2\pi} \right)^{1-\sigma-i\xi} 
\nonumber
\end{equation}
\begin{equation}
+    
\frac{i}{\pi} \left( \frac{\eta}{2\pi} \right)^{-(\sigma+i\xi)} \Biggl( -i \arg{(1-e^{i\eta})} + \frac{\xi -i\sigma}{\eta} \Re{\operatorname{Li_2}{\left(e^{i\eta}\right)}} 
\nonumber
\end{equation}
\begin{equation}
+\frac{1}{\eta^2} [i\xi^2 + (2\sigma+1)\xi -i\sigma(\sigma+1)] \Im{\operatorname{Li_3}{\left(e^{i\eta}\right)}}  \Biggr) 
\nonumber
\end{equation}
\begin{equation}
+ 
\left( \frac{\eta}{2\pi} \right)^{-(\sigma+i\xi)} O\left( \frac{\xi^3}{\eta^{3+\sigma}} \right), \quad \eta>\xi, \quad 0\le\sigma\le 1, \quad \xi\to\infty,
\end{equation}
where $\operatorname{Li_m}(z)$ denotes the polylogarithm, i.e.,
\begin{equation}
\operatorname{Li_m}(e^{i\eta}) = \sum_{k=1}^{\infty} \frac{e^{ik\eta}}{k^m}, \quad m\ge 1.
\end{equation}
It should be emphasized that the {\it only} oscillatory dependence of $\xi$ in the last term in the rhs of \eqref{5.10} is in the form $(\eta/2\pi)^{-(\sigma +i\xi)}$.

Letting $\xi=t-x$ and $\eta=t$ in \eqref{5.10} we find
\begin{multline*}
\zeta(\sigma +it - ix) = \sum_{m=1}^{\left \lfloor{T}\right \rfloor} \frac{1}{m^{\sigma+it-ix}} - \frac{1}{1 - \sigma - it + ix} ~ T^{1 - \sigma - it + ix} \\
+ \frac{i}{\pi} ~ T^{- \sigma - it + ix} \Big( -i \arg (1 - e^{it}) + \frac{t - x - i\sigma}{t} \Re{\operatorname{Li_2}{(e^{it})}} \\
+ \frac{1}{t^2} \left[  it^2 + ix^2 -2itx + (2\sigma + 1) t - (2\sigma + 1) x -i\sigma (\sigma + 1)  \right] \Im{\operatorname{Li_3}{(e^{it})}}   \Big) \\
+  T^{- \sigma - it + ix} O \left(  \frac{(t-x)^3}{t^{3+\sigma}} \right), \qquad |x| \le t^{\delta_{34}}, \qquad t\to\infty.
\end{multline*}

Hence, using the fact that
\begin{equation*}
\operatorname{Li_m}{(e^{it})} = O(1), \quad t\to\infty, \quad m\ge 1,
\end{equation*}
we obtain
\begin{subequations} \label{5.12}
\begin{equation}  \label{5.12a}
\zeta(\sigma +it - ix) = \sum_{m=1}^{\left \lfloor{T}\right \rfloor} \frac{m^{ix}}{m^s} + G(e^{it}) \frac{T^{- \sigma - it + ix}}{\pi} + O \left( t^{-\sigma}, \frac{x}{t}  \right) T^{- \sigma - it + ix},
\end{equation}
where $G(e^{it})$ is defined by
\begin{equation}  \label{5.12b}
G(e^{it}) = \arg{(1-e^{it})} - \Im{\operatorname{Li_3}{(e^{it})}} - \frac{i}{2} + i \Re{\operatorname{Li_2}{(e^{it})}}.
\end{equation}
\end{subequations}
Hence,
\begin{multline} \label{5.13}
|\zeta(\sigma +it - ix) |^2 =  \sum_{m_1=1}^{\left \lfloor{T}\right \rfloor}  \sum_{m_2=1}^{\left \lfloor{T}\right \rfloor}  \frac{1}{m_1^s m_2^{\bar{s}}} \left(\frac{m_1}{m_2}\right)^{ix} + \sum_{m=1}^{T} \frac{m^{ix}}{m^s}T^{-ix} \overline{G(e^{it})} T^{it} \frac{T^{-\sigma}}{\pi} \\
\qquad + \sum_{m=1}^{\left \lfloor{T}\right \rfloor} \frac{m^{-ix}}{m^{\bar{s}}}T^{ix} G(e^{it}) T^{-it} \frac{T^{-\sigma}}{\pi} + |G|^2 \frac{T^{-2\sigma}}{\pi^2} + \sum_{m=1}^{\left \lfloor{T}\right \rfloor} \frac{m^{ix}}{m^{s}}T^{-ix} O \left( t^{-2\sigma}, \frac{x}{t} t^{-\sigma}  \right) T^{it} \\
+ \left(  G(e^{it}) + \overline{G(e^{it})} \right) T^{-2\sigma} O \left(  t^{-\sigma}, \frac{x}{t} \right) + T^{-2\sigma} O \left(  t^{-2\sigma}, \frac{x}{t}t^{-\sigma}, \left( \frac{x}{t}  \right)^2 \right)\\
+ \sum_{m=1}^{\left \lfloor{T}\right \rfloor} \frac{m^{-ix}}{m^{\bar{s}}}T^{ix} O \left( t^{-2\sigma}, \frac{x}{t} t^{-\sigma}  \right) T^{-it}, \quad |x|\le t^{\delta_{34}}, \quad t\to\infty.
\end{multline}
Substituting \eqref{5.9} and \eqref{5.13} into \eqref{5.6} we find
\begin{align} \label{5.14}
I_4(\sigma, t,\delta_3,\delta_4) &= I_4^{(1)}(\sigma, t,\delta_3,\delta_4) + \overline{G(e^{it})}T^{it} \frac{T^{-\sigma}}{\pi} I_4^{(2)}(\sigma, t,\delta_3,\delta_4) \nonumber \\
& \quad+ G(e^{it}) T^{-it} \frac{T^{-\sigma}}{\pi} I_4^{(3)}(\sigma, t,\delta_3,\delta_4) + |G|^2 \frac{T^{-2\sigma}}{\pi^2} I_4^{(4)}(\sigma, t,\delta_3,\delta_4) \nonumber \\
& \quad+  I_4^{(5)}(\sigma, t,\delta_3,\delta_4) + I_4^{(6)}(\sigma, t,\delta_3,\delta_4),
\end{align}
where
\begin{subequations}  \label{5.15}
\begin{multline}  \label{5.15a}
I_4^{(1)}(\sigma, t,\delta_3,\delta_4) = \frac{1}{\pi} \Re{\left\{ \sum_{m_1=1}^{\left \lfloor{T}\right \rfloor}  \sum_{m_2=1}^{\left \lfloor{T}\right \rfloor} \frac{1}{m_1^s m_2^{\bar{s}}} \int_{H} \frac{e^z}{z}J_4(\sigma, t,\delta_3,\delta_4,A_1)\textrm{d}z \right\}} \\
\times \left[ 1 + O\left( \frac{1}{t} \right) \right], \quad t\to\infty,
\end{multline}
\begin{multline} \label{5.15b}
I_4^{(2)}(\sigma, t,\delta_3,\delta_4) = \frac{1}{\pi} \Re{\left\{\sum_{m=1}^{\left \lfloor{T}\right \rfloor} \frac{1}{m^s} \int_{H} \frac{e^z}{z}J_4(\sigma,t,\delta_3,\delta_4,A_2)\textrm{d}z \right\}} \\
\times \left[ 1 + O\left( \frac{1}{t} \right) \right], \quad t\to\infty,
\end{multline}
\begin{multline} \label{5.15c}
I_4^{(3)}(\sigma, t,\delta_3,\delta_4) = \frac{1}{\pi} \Re{\left\{\sum_{m=1}^{\left \lfloor{T}\right \rfloor} \frac{1}{m^{\bar{s}}} \int_{H} \frac{e^z}{z}J_4(\sigma,t,\delta_3,\delta_4,A_3)\textrm{d}z \right\}} \\
\times \left[ 1 + O\left( \frac{1}{t} \right) \right], \quad t\to\infty,
\end{multline}
\begin{multline} \label{5.15d}
I_4^{(4)}(\sigma, t,\delta_3,\delta_4) = \frac{1}{\pi} \Re{\left\{ \int_{H} \frac{e^z}{z}J_4(\sigma,t,\delta_3,\delta_4,A_4)\textrm{d}z \right\}} \\
\times \left[ 1 + O\left( \frac{1}{t} \right) \right], \quad t\to\infty,
\end{multline}
\begin{equation} \label{5.15e}
I_4^{(5)}(\sigma, t,\delta_3,\delta_4) = \frac{1}{\pi} \Re{\left\{\sum_{m=1}^{\left \lfloor{T}\right \rfloor} \frac{1}{m^s} \int_{H} \frac{e^z}{z}\hat{J}_4(\sigma,t,\delta_3,\delta_4,A_5)\textrm{d}z \right\}},
\end{equation}
\begin{equation} \label{5.15f}
I_4^{(6)}(\sigma, t,\delta_3,\delta_4) = \frac{1}{\pi} \Re{\left\{\sum_{m=1}^{\left \lfloor{T}\right \rfloor} \frac{1}{m^{\bar{s}}} \int_{H} \frac{e^z}{z}\hat{J}_4(\sigma,t,\delta_3,\delta_4,A_6)\textrm{d}z \right\}},
\end{equation}
\end{subequations}
where
\begin{equation} \label{5.17}
J_4(\sigma,t,\delta_3,\delta_4,A) =  
\nonumber
\end{equation}
\begin{multline}
PV \int_{-t^{\delta_3}}^{t^{\delta_4}} \frac{e^{\frac{\pi x}{2}}A^{ix}}{e^{-\pi x}-e^{\pi x}} e^{ix} \left(  1 - \frac{x}{t} \right)^{i(t-x)} \left[ 1 + O\left( \frac{1}{t-x} \right) \right] \left( 1 - \frac{x}{t} \right)^{\sigma-\frac{1}{2}} \textrm{d}x, \\
t\to\infty,
\end{multline}
\begin{multline}
\hat{J}_4(\sigma,t,\delta_3,\delta_4,A) = PV \int_{-t^{\delta_3}}^{t^{\delta_4}} \frac{e^{\frac{\pi x}{2}}}{e^{-\pi x} - e^{\pi x}} A^{ix} e^{ix} \left( 1 - \frac{x}{t}  \right)^{i(t-x)} \\
\times O\left(  t^{-2\sigma}, \frac{x}{t} t^{-\sigma} \right) \left[ 1 + O \left(  \frac{1}{t-x} \right)  \right] \left(  1 - \frac{x}{t} \right)^{\sigma - \frac{1}{2}}\textrm{d}x, \quad t\to\infty,
\end{multline}
with
\begin{equation} \label{5.18}
A_1=\frac{m_1}{m_2} \frac{z}{t}, \quad A_2=2\pi m\frac{z}{t^2}, \quad A_3 = \frac{z}{2\pi m}, \quad A_4 = \frac{z}{t}, \quad A_5 = A_2, \quad A_6=A_3.
\end{equation}
Using the identity
\begin{align*} \label{5.19}
e^{ix}\left(  1 - \frac{x}{t} \right)^{i(t-x)} &= e^{ix} e^{i(t-x)\ln{\left(  1 - \frac{x}{t} \right)}} =e^{ix} e^{i(t-x) \left(  -\frac{x}{t} + O \left(  \frac{x^2}{t^2} \right) \right) } \\
&= e^{iO \left(  \frac{x^2}{t} \right) } = 1 + O \left(  \frac{x^2}{t} \right), \quad \frac{x}{t} \to 0,
\end{align*}
we find
\begin{multline}
e^{ix} \left(  1 - \frac{x}{t} \right)^{i(t-x)} \left[ 1 + O\left( \frac{1}{t-x} \right) \right] \left( 1 - \frac{x}{t} \right)^{\sigma - \frac{1}{2}} \\
= \left[ 1 + O\left( \frac{t^{\delta_{34}}}{t} \right) \right], \quad |x|<t^{\delta_{34}}, \quad t\to\infty.
\end{multline}
Thus, equation \eqref{5.17} yields
\begin{equation} \label{5.20}
J_4(\sigma, t,\delta_3,\delta_4,A) = \tilde{J_4}(\sigma, t,\delta_3,\delta_4,A) \left[  1 + O \left( \frac{t^{\delta_{34}}}{t} \right) \right], \quad t\to\infty,
\end{equation}
where $\tilde{J_4}$ is defined by
\begin{equation} \label{5.21}
\tilde{J_4}(\sigma, t,\delta_3,\delta_4,A) = PV\int_{-t^{\delta_3}}^{t^{\delta_4}} \frac{e^{\frac{\pi x}{2}}A^{ix}}{e^{-\pi x}-e^{\pi x}} \textrm{d}x.
\end{equation}

It is shown in proposition 5.1 of \cite{F} that
\begin{equation} \label{5.30}
\tilde{J_4}(t,\delta_3,\delta_4,A)  = \frac{i}{2} \left( -1 + \frac{2}{1-iA} \right) + O\left( e^{-t^{\tilde{\delta}_{34}}} \right), \quad t \to\infty.
\end{equation}
In order to compute the integral $\tilde{J_4}$ around the Hankel contour occurring in \eqref{5.15} we will employ the following residue formulae:
\begin{equation} \label{5.31}
\int_{H} \frac{e^z}{z}\textrm{d}z = 2\pi i,
\end{equation}
and
\begin{equation} \label{5.32}
\int_{H} \frac{e^z}{z(z+ic)}\textrm{d}z = 2\pi i  \underset{z=0}{\operatorname{Res}}{\frac{e^z}{z(z+ic)}} = \frac{2\pi}{c}, \quad c\ne 0.
\end{equation}
These formulae imply
\begin{equation} \label{5.33}
-\frac{i}{2}\int_{H} \frac{e^z}{z}\textrm{d}z = \pi,
\end{equation}
as well as
\begin{equation} \label{5.34}
i\int_{H} \frac{e^z}{z(1-iA_1)}\textrm{d}z = \frac{i}{-i\left( \frac{m_1}{m_2t} \right)} \int_{H} \frac{e^z}{z\left(z+ \frac{im_2t}{m_1} \right)}\textrm{d}z = -2\pi,
\end{equation}
and
\begin{equation} \label{5.35}
i\int_{H} \frac{e^z}{z(1-iA_2)}\textrm{d}z = -2\pi.
\end{equation}
Hence,
\begin{equation}\label{5.36}
\int_{H} \frac{e^z}{z} \tilde{J_4}(t,\delta_3,\delta_4,A_j) \textrm{d}\tau = -\pi + O\left( e^{-t^{\tilde{\delta}_{34}}} \right), \quad j=1,2.
\end{equation}

Similar considerations apply to the calculation of $\tilde{J}_4(t, \delta_3, \delta_4, A_j)$, $j=$ 3, 4, 5, 6.

Using \eqref{5.36} in \eqref{5.15} and then substituting the resulting expressions in \eqref{5.14} we find
\begin{equation}  \label{5.3}
I_4(\sigma, t,\delta_3,\delta_4) = - \sum_{m_1=1}^{\left \lfloor{T}\right \rfloor}  \sum_{m_2=1}^{\left \lfloor{T}\right \rfloor} \frac{1}{m_1^s m_2^{\bar{s}}}  \left[ 1 + O\left( \frac{t^{\delta_{34}}}{t} \right) \right]
\nonumber
\end{equation}
\begin{equation}
+ \Re{\left\{ \sum_{m=1}^{\left \lfloor{T}\right \rfloor} \frac{1}{m^s} \right\}}  O\left( t^{-\sigma} \right), \quad t\to\infty.
\end{equation}
In order to estimate the single and double sums appearing in \eqref{5.3} we will use the following ``crude" estimates:
\begin{align}\label{5.40}
  \Bigg| \sum_{m_1=1}^{[T]}  \sum_{m_2=1}^{[T]}\frac{1}{m_1^s m_2^{\bar{s}}}  \Bigg| & \leq \int_1^T \int_1^T\dfrac{1}{x^\sigma} \dfrac{1}{y^\sigma} dxdy =O\left(t^{2-2\sigma}\right), \nonumber \\  \Bigg| \sum_{m=1}^{[T]}\frac{1}{m^s} \Bigg| & \leq \int_1^T \dfrac{1}{x^\sigma} dx =O\left(t^{1-\sigma}\right).
\end{align}  
Then, equation \eqref{5.3} yields \eqref{5.3n}.

\textbf{QED}

\vphantom{a}

\begin{remark}\label{rem5}

By employing techniques developed in \cite{T} and \cite{T2} it is possible to improve the estimates in \eqref{5.40} to estimates of $\displaystyle O\left(t^{\frac{1}{2}-\frac{3}{2}\sigma+\delta_{34}}\right),$ see \cite{FK}.

\end{remark}

\section{The Analysis of $I_{3}$} \label{sec6}

\begin{lemma}
Let $I_{3}(\sigma ,t,\delta _{2},\delta _{3})$ be defined by

\begin{equation} \label{sixone}
I_{3}(\sigma ,t,\delta _{2},\delta _{3})=\frac{t}{\pi}
\int_{t^{   \delta _{2} -1  }}^{1-t^{\delta _{3}-1}}
\Re\left\{ \frac{\Gamma(it - i\tau t) }{\Gamma(\sigma + i t)}  \Gamma(\sigma + i\tau t) \right\}
|\zeta (\sigma +it\tau)|^{2}d\tau , \nonumber 
\end{equation}
\begin{equation}
0<\sigma <1, \hphantom{2a} t>0,
\end{equation}

where $\zeta (z)$ and $\Gamma (z)$ denote the Riemann zeta and gamma functions respectively, and $\delta _{2}$, $\delta _{3}$ 
are sufficiently small, positive constants. Then,
\begin{equation}
I_{3}(\sigma ,t,\delta _{2},\delta _{3})=\sqrt{\frac{2t}{\pi}}
\Re\Biggl(  
e^{-\frac{i\pi}{4}}\sum_{m_{1}=1}^{[T]}\sum_{m_{2}=1}^{[T]}m_{1}^{-\sigma }m_{2}^{-\sigma }
J^{(1)}(\sigma ,t,\delta _{2},\delta _{3},\frac{m_{2}}{m_{1}}) \nonumber
\end{equation}
\begin{equation}
+T^{-\sigma }e^{-\frac{i\pi}{4}}\sum_{m=1}^{[T]}m^{-\sigma }[
\frac{i}{2\pi}J^{(2)}(\sigma ,t,\delta _{2},\delta _{3},\frac{T}{m})
-\frac{i}{2\pi}J^{(3)}(\sigma ,t,\delta _{2},\delta _{3},\frac{m}{T})  \nonumber
\end{equation}
\begin{equation}
+J^{(4)}(\sigma ,t,\delta _{2},\delta _{3},\frac{T}{m})
+J^{(5)}(\sigma ,t,\delta _{2},\delta _{3},\frac{m}{T})
]+T^{-2\sigma }e^{-\frac{i\pi}{4}}J^{(6)}(\sigma ,t,\delta _{2},\delta _{3})\Biggr)  \nonumber
\end{equation}
\begin{equation}  \label{sixtwo}
\times [1+O(\frac{1}{t})],   \hphantom{2a}
0<\sigma <1, \hphantom{2a} t\rightarrow \infty ,
\end{equation}
where $\{ J^{(j)} \}_1^6$ are defined as follows:
\begin{equation} \label{sixthree}
J^{(1)}(\sigma ,t,\delta _{2},\delta _{3},\lambda )=
\int_{t^{   \delta _{2} -1  }}^{1-t^{\delta _{3}-1}}
G(\sigma ,\tau )e^{itF(\tau ,\lambda )}A(t,\tau )d\tau ,
\end{equation}
\begin{equation} \label{sixfour}
J^{(2)}(\sigma ,t,\delta _{2},\delta _{3},\lambda )=
\int_{t^{   \delta _{2} -1  }}^{1-t^{\delta _{3}-1}}
G(\sigma ,\tau )e^{itF(\tau ,\lambda )}\frac{A(t,\tau )}{\tau - i\frac{(1-\sigma )}{t}}d\tau ,
\end{equation}
\begin{equation} \label{sixfive}
J^{(3)}(\sigma ,t,\delta _{2},\delta _{3},\lambda )=
\int_{t^{   \delta _{2} -1  }}^{1-t^{\delta _{3}-1}}
G(\sigma ,\tau )e^{itF(\tau ,\lambda )}\frac{A(t,\tau )}{\tau + i\frac{(1-\sigma )}{t}}d\tau ,
\end{equation}
\begin{equation} \label{sixsix}
J^{(4)}(\sigma ,t,\delta _{2},\delta _{3},\lambda )=
\int_{t^{   \delta _{2} -1  }}^{1-t^{\delta _{3}-1}}
G(\sigma ,\tau )e^{itF(\tau ,\lambda )}A(t,\tau )\overline{B(\sigma,t,\tau )}d\tau ,
\end{equation}
\begin{equation} \label{sixseven}
J^{(5)}(\sigma ,t,\delta _{2},\delta _{3},\lambda )=
\int_{t^{   \delta _{2} -1  }}^{1-t^{\delta _{3}-1}}
G(\sigma ,\tau )e^{itF(\tau ,\lambda )}A(t,\tau )B(\sigma, t,\tau )d\tau,
\end{equation}
\begin{equation} \label{sixeight}
J^{(6)}(\sigma ,t,\delta _{2},\delta _{3} )=
\int_{t^{   \delta _{2} -1  }}^{1-t^{\delta _{3}-1}}
G(\sigma ,\tau )e^{itF(\tau ,1 )}A(t,\tau )C(\sigma,t,\tau )d\tau,
\end{equation}
with
\begin{equation}  \label{sixnine}
G(\sigma ,\tau )=(1-\tau )^{-\frac{1}{2}}\tau ^{\sigma -\frac{1}{2}},
\end{equation}
\begin{equation}  \label{sixten}
F(\tau ,\lambda )=(1-\tau )\ln (1-\tau )+\tau \ln \tau +\tau \ln\lambda ,
\end{equation}
\begin{equation}  \label{sixeleven}
A(t,\tau )=1+O(\frac{1}{t-t\tau })+O(\frac{1}{t\tau }),  \hphantom{2a} t\rightarrow \infty ,
\end{equation}
\begin{equation} \label{sixtwelve}
B(\sigma,t,\tau )=O(1,\tau ,\tau ^{2},\frac{\tau ^{3}}{t^{\sigma }}),  \hphantom{2a} t\rightarrow \infty ,
\end{equation}
\begin{equation}  \label{sixthirteen}
C(\sigma,t,\tau )=\frac{1}{(2\pi)^{2}}\frac{1}{\tau ^{2}+(\frac{1-\sigma }{t})^{2}}
+O(\frac{B}{\tau \pm i\frac{1-\sigma }{t}},|B|^{2}),  \hphantom{2a} t\rightarrow \infty .
\end{equation}
\end{lemma}

\textbf{Proof.} In the interval of integration we have that
\begin{equation}
t^{\delta _{2}}\leq t\tau \leq t-t^{\delta _{3}}. \nonumber
\end{equation}
Thus, $t-t\tau\to\infty$ as $t\to\infty$ and hence we can use the following asymptotic expansion for $\Gamma (it -it\tau )/\Gamma (\sigma +it )$:

\begin{align} \label{4.9}
\frac{\Gamma(it -it\tau)}{\Gamma(\sigma + i t)} &= \frac{[(1-\tau)t]^{i(1-\tau)t}[(1-\tau)t]^{-\frac{1}{2}}e^{-i(1-\tau)t}e^{-\frac{\pi}{2}(1-\tau)t}\left[ 1 + O\left( \frac{1}{t-t\tau} \right) \right]} {t^{it}t^{\sigma-\frac{1}{2}}e^{-it}e^{\frac{i\pi\sigma}{2}}e^{-\frac{\pi t}{2}}\left[ 1 + O\left( \frac{1}{t} \right) \right]} \nonumber \\
&=t^{-\sigma}e^{\frac{-i\pi\sigma}{2}}(1-\tau)^{i(1-\tau)t} t^{-i\tau t} (1-\tau)^{-\frac{1}{2}} e^{i\tau t}e^{\frac{\pi\tau t}{2}}
\nonumber
\end{align}
\begin{equation} \label{49}
\times \left[ 1 + O\left( \frac{1}{t-t\tau} \right) \right] \left[ 1 + O\left( \frac{1}{t} \right) \right], ~ t\to\infty.
\end{equation}
Furthermore, since 
\begin{equation}  \label{sixfourteen}
t\tau \rightarrow \infty \hphantom{3a} \text{as} \hphantom{3a} t \rightarrow \infty, 
\end{equation}
we can employ the asymptotic formulae \eqref{219}  with $\xi =t\tau $ to evaluate $\Gamma (\sigma +it\tau )$:
\begin{equation}  \label{sixfifteen}
\Gamma (\sigma +it\tau )=\sqrt{2\pi}e^{-\frac{i\pi}{4}}(t\tau )^{\sigma -\frac{1}{2}}e^{-\frac{\pi t\tau }{2}}
(t\tau )^{it\tau }e^{-it\tau }e^{\frac{i\pi\sigma }{2}}[1+O(\frac{1}{t\tau })], \hphantom{2a} t\tau \rightarrow \infty .
\end{equation}
Equations \eqref{49} and (\ref{sixfifteen}) imply
\begin{equation} 
\frac{\Gamma (it-it\tau )}{\Gamma (\sigma +it )}\Gamma (\sigma +it\tau )=
\sqrt{\frac{2\pi}{t}}e^{-\frac{i\pi}{4}}G(\sigma ,\tau )e^{it[(1-\tau )\ln (1-\tau )+\tau \ln \tau  ]}A(t,\tau )
[1+O(\frac{1}{t})],   \nonumber
\end{equation}
\begin{equation}  \label{sixsixteen}
 \hphantom{2a} t-t\tau \rightarrow \infty ,  \hphantom{2a} \tau t \rightarrow \infty 
, \hphantom{2a}   t \rightarrow \infty ,  
\end{equation}
where $G(\sigma ,\tau )$ and $A(t,\tau )$ are defined in (\ref{sixnine}) and (\ref{sixeleven}). 

Equation (\ref{sixfourteen}) implies that we can employ equation \eqref{5.10} with $\eta =t$ and $\xi =t\tau $:
\begin{equation}
\zeta (\sigma +it\tau )=\sum_{m=1}^{[T]}\frac{1}{m^{\sigma +it\tau }} - \frac{i}{2\pi}
\frac{T^{-(\sigma +it\tau )}}{\tau + i\frac{(1-\sigma )}{t}}+\frac{i}{\pi}T^{-(\sigma +it\tau )}\Biggl( -i\arg (1-e^{it})
\nonumber
\end{equation}
\begin{equation}
+(\tau -i\frac{\sigma }{t})\Re Li_{2}(e^{it})+[i\tau ^{2}+(2\sigma +1)\frac{\tau }{t}-i\frac{\sigma (\sigma +1)}{t^{2}}]
\Im Li_{3}(e^{it}) \Biggr)  \nonumber
\end{equation}
\begin{equation}  
+ T^{-(\sigma +it\tau )}O(\frac{\tau ^{3}}{t^{\sigma }}), \hphantom{3a} 
\frac{t^{\delta _{2}}}{t}\leq\tau \leq 1-\frac{t^{\delta _{3}}}{t}, \hphantom{2a} t\rightarrow \infty . \nonumber 
\end{equation}
Hence,
\begin{subequations}
\begin{equation}  \label{sixseventeen}
\zeta (\sigma +it\tau )=\sum_{m=1}^{[T]}\frac{1}{m^{\sigma +it\tau }} - \frac{i}{2\pi}
\frac{T^{-(\sigma +it\tau )}}{\tau + i\frac{(1-\sigma )}{t}}+T^{-(\sigma +it\tau )}B(t,\tau )
, \hphantom{2a} t\rightarrow \infty ,
\end{equation}
where $B$ is defined by
\begin{equation}
B(\sigma,t,\tau )= \nonumber
\end{equation}
\begin{equation}  
\frac{i}{\pi}\Biggl( -i\arg (1-e^{it})
+(\tau -i\frac{\sigma }{t})\Re Li_{2}(e^{it})+[i\tau ^{2}+(2\sigma +1)\frac{\tau }{t}-i\frac{\sigma (\sigma +1)}{t^{2}}]
\Im Li_{3}(e^{it}) \Biggr)+  \nonumber
\end{equation}
\begin{equation}   \label{sixseighteen}
O(\frac{\tau ^{3}}{t^{\sigma }}).
\end{equation}
\end{subequations}

Multiplying the expression (\ref{sixseventeen}) with the expression obtained from (\ref{sixseventeen}) via complex conjugation, 
i.e., with the expression 
\begin{equation} 
\overline{ \zeta (\sigma +it\tau )}=
\sum_{m=1}^{[T]}\frac{1}{m^{\sigma -it\tau }} + \frac{i}{2\pi}
\frac{T^{-\sigma +it\tau }}{\tau - i\frac{(1-\sigma )}{t}}+T^{-\sigma +it\tau }\overline{B(t,\tau )},  \nonumber
\end{equation}
we find
\begin{equation}  \label{sixnineteen}
|\zeta (\sigma +it\tau )|^{2}=
\sum_{m_{1}=1}^{[T]}\sum_{m_{2}=1}^{[T]}
\frac{1}{m_{1}^{\sigma }m_{2}^{\sigma }}(\frac{m_{2}}{m_{1}})^{it\tau }+ \nonumber
\end{equation}
\begin{equation}
T^{-\sigma }\sum_{m=1}^{[T]}\frac{1}{m^{\sigma }}\left\{ 
\frac{i}{2\pi}\frac{(\frac{T}{m})^{it\tau }}{\tau - i\frac{(1-\sigma )}{t}}
- \frac{i}{2\pi}\frac{(\frac{m}{T})^{it\tau }}{\tau + i\frac{(1-\sigma )}{t}}
+(\frac{T}{m})^{it\tau }\bar{B}+(\frac{m}{T})^{it\tau }B 
 \right\}+  \nonumber
\end{equation}
\begin{equation}
T^{-2\sigma }C(\sigma,t,\tau ),  \hphantom{2a} t\rightarrow \infty ,
\end{equation}
where $C$ is defined by
\begin{equation}  \label{sixtwenty}
C(\sigma, t,\tau )=\frac{1}{(2\pi)^{2}}\frac{1}{\tau ^{2}+\frac{(1-\sigma )^{2}}{t^{2}}}+|B|^{2}
+\frac{i}{2\pi}\frac{B}{\tau - i\frac{(1-\sigma )}{t}}
- \frac{i}{2\pi}\frac{\bar{B}}{\tau + i\frac{(1-\sigma )}{t}}.
\end{equation}
We note that
\begin{equation}
B(\sigma,t,\tau )=O(1,\tau ,\tau ^{2},\frac{\tau }{t},\frac{\tau ^{3}}{t^{\sigma }}),  \hphantom{2a} t\rightarrow \infty , 
\nonumber
\end{equation}
whereas
\begin{equation}
C(\sigma, t,\tau )=\frac{1}{(2\pi)^{2}}\frac{1}{\tau ^{2}+\frac{(1-\sigma )^{2}}{t^{2}}}+
O(\frac{B}{\tau \pm i\frac{1-\sigma }{t}},|B|^{2}),  \hphantom{2a} t\rightarrow \infty .  \nonumber
\end{equation}
Substituting equations (\ref{sixsixteen}) and (\ref{sixnineteen}) in equation (\ref{sixone}) we find equation (\ref{sixtwo}). 

\textbf{QED}

\vphantom{a}

It is well known that the main contributions of the asymptotic analysis of integrals come from possible singularities, from possible stationary points, and from the end points of the interval of integration \cite{AF}. Each of the integrals $\left\{J^{(j)}\right\}_1^5$ possesses a stationary point at $\tau =1/(1+\lambda )$, whereas for $J^{(6)}$ the possible stationary point occurs at $\tau =1$ which is outside the interval of integration. Thus, for the integrals $\left\{ J^{(j)} \right\}_1^5$ there exist two contributions, one from the associated stationary points and one from the end points of the interval of integration, whereas for $J^{(6)}$ there exists only the latter contribution.

The equation
\begin{equation*}
G \left( \sigma, 1 - \frac{t^{\delta_3}}{t}  \right) = t^{\frac{1 - \delta_3}{2}} \left( 1 - \frac{t^{\delta_3}}{t}  \right)^{\sigma - \frac{1}{2}},
\end{equation*}
implies that there exists an integrable singularity in the neighborhood of the end point $1-t^{\delta_3}/t$, thus the associated contribution is larger than the usual contribution of order $O(1/t)$ as $t\to\infty$.

The equation
\begin{equation*}
G \left( \sigma, \frac{t^{\delta_2}}{t}  \right) = t^{\left(1-\delta_2\right)\left(\frac{1}{2}-\sigma\right)} \left( 1 - \frac{t^{\delta_2}}{t}  \right),
\end{equation*}
implies that if $\sigma\ge 1/2$, the contribution of the end point $t^{\delta_2}/t$ can be computed following the standard integration by parts arguments, and it is of order $O(1/t)$ as $t\to\infty$. Also, for $\sigma>0$, the contribution of the end point $t^{\delta_2}/t$ is always smaller than the contribution of the end point $1-t^{\delta_3}/t$. However, for economy of presentation we will assume that $\sigma\ge 1/2$, thus we only have to compute the contributions of the stationary points of the integrals $\left\{J^{(j)}\right\}_1^5$, as well as the contribution from the end point of all integrals $\left\{J^{(j)}\right\}_1^6$.

In Lemma \ref{l6.2a}, we compute the former contribution; the latter contribution will be computed in Lemma \ref{l6.2}.

\begin{lemma}\label{l6.2a}
Let the integrals $\{ J^{(j)}(\sigma ,t,\delta _{2},\delta _{3}, \lambda ) \}_{1}^{5}$ be defined in 
(\ref{sixthree})-(\ref{sixseven}). Then,

\begin{equation} 
J^{(j)}(\sigma ,t,\delta _{2},\delta _{3}, \lambda )=
\sqrt{\frac{2\pi}{t}}e^{\frac{i\pi}{4}}\frac{\lambda ^{it}}{(1+\lambda )^{\sigma +it}}P_{j}(\lambda ,t)
\left[1+O\left(\frac{1+\frac{1}{\lambda} }{t}\right)+O\left(\frac{1+\lambda}{t}\right)\right]  \nonumber 
\end{equation}
\begin{equation}  \label{sixtwentyone}
+J_{B}^{(j)}(\sigma ,t,\delta _{3}, \lambda ), \quad j=1,\dots ,5, \quad \frac{1}{2}\le\sigma <1, \quad t\rightarrow \infty , 
\end{equation}
where the integrals $\{ J_{B}^{(j)} \}_{1}^{5}$ are obtained from the integrals $\{ J^{(j)} \}_{1}^{5}$ 
defined in (\ref{sixthree})-(\ref{sixseven}) by replacing the contour of integration with the ray from the point 
$1-t^{\delta _{3}}/t$ to $\infty e^{i\varphi }$ with $\varphi $ satisfying 
\begin{equation} \label{sixtwentytwo}
0<\varphi <\arctan (\frac{\pi}{|\ln\lambda |}),
\end{equation}
and $ \left\{ P_{j}(\lambda ,t) \right\}_1^5$ defined as follows:
\begin{equation} \label{sixtwentythree}
P_{1}(\lambda ,t)=1,
\end{equation}
\begin{equation} \label{sixtwentyfour}
P_{2}(\lambda ,t)=\frac{1+\lambda }{1-\frac{i(1-\sigma )}{t}(1+\lambda )},
\end{equation}
\begin{equation} \label{sixtwentyfive}
P_{3}(\lambda ,t)=\frac{1+\lambda }{1+\frac{i(1-\sigma )}{t}(1+\lambda )},
\end{equation}
\begin{equation} \label{sixtwentysix}
P_{4}(\lambda ,t)=O\left(1,\frac{1}{1+\lambda },\frac{1}{(1+\lambda )^{2}},\frac{1}{(1+\lambda )^{3}t^{\sigma }}\right), \quad t\to \infty,
\end{equation}
\begin{equation} \label{sixtwentyseven}
P_{5}(\lambda ,t)=P_{4}(\lambda ,t), \quad t\to \infty.
\end{equation}
In the above formulae $\lambda $ satisfies
\begin{equation}   \label{sixtwentyeight}
\frac{1}{t^{1-\delta _{3}}-1}<\lambda <t^{1-\delta _{2}}-1.
\end{equation}
\end{lemma}

\textbf{Proof.} The definition of $F$ given in (\ref{sixten}) implies 
\begin{equation} \label{sixtwentynine}
\frac{\partial F}{\partial \tau }=-\ln (1-\tau )+\ln \tau +\ln \lambda .
\end{equation}
Thus, a possible stationary point occurs at $\tau =\tau _{1}$, where 
\begin{equation} \label{sixtwenthirty}
\tau _{1}=\frac{1}{1+\lambda }.
\end{equation}
Thus, s stationary point occurs {\it inside} the interval of integration if and only if
\begin{equation*} 
\frac{t^{\delta _{2}}}{t}<\tau _{1}<1-\frac{t^{\delta _{3}}}{t},
\end{equation*}
i.e. if and only if $\lambda $ satisfies the inequality (\ref{sixtwentyeight}).

We deform the contour of integration of $\{ J^{(j)} \}_{1}^{5}$, from the point $t^{\delta _{2}-1}$ down into the lower half of the 
complex $\tau $-plane, then up through the point $\tau =\tau _{1}$ continuing to $\infty e^{i\varphi }$, $\varphi >0$, and 
finally back to the point $1-t^{\delta _{3}-1}$.

We claim that if $\varphi $ is sufficiently small, namely if $\varphi $ satisfies (\ref{sixtwentytwo}), then the integrals 
$\{ J_{B}^{(j)} \}_{1}^{5}$ converge. Indeed, employing the change of variables
\begin{equation*}
\tau =\Delta (t)+\rho e^{i\varphi }, \hphantom{3a} \Delta(t) =1-\frac{t^{\delta _{3}}}{t},
\end{equation*}
we find that $F$ becomes 
\begin{equation*}
F = (1-\Delta -e^{i\varphi }\rho )\ln (1-\Delta -e^{i\varphi }\rho ) + (\Delta +e^{i\varphi }\rho )[\ln\lambda +\ln (\Delta +e^{i\varphi }\rho )].
\end{equation*}
For fixed $t$,
\begin{equation}
F\sim \rho e^{i\varphi }[\ln \lambda +\ln (\rho e^{i\varphi })-\ln (-\rho e^{i\varphi })], \hphantom{3a} 
\rho \rightarrow \infty .   \nonumber
\end{equation}
Using that
\begin{equation*}
\ln{(\rho e^{i\varphi})} - \ln{(-\rho e^{i\varphi})} = i\pi,
\end{equation*}
it follows that
\begin{equation*}
\Im F\sim \rho [\ln\lambda \sin\varphi +\pi\cos\varphi ],  \hphantom{3a} \rho \rightarrow \infty .
\end{equation*}
For the convergence of $J_{B}^{(j)}$ we require $\Im F>0$ as $\rho \rightarrow \infty $. Thus, if $\lambda \geq 1$ then we have convergence for all $\varphi \in (0,\pi/2)$, whereas if $\lambda \in (0,1)$ we require that $\varphi $ satisfies 
(\ref{sixtwentytwo}).

In order to compute the contribution from the stationary point $\tau =\tau _{1}$ we employ the well known formula \cite{M}
\begin{equation}  \label{sixthirtyone}
\int g(\tau )e^{itF(\tau )}d\tau =
\sqrt{\frac{2\pi}{t|F^{\prime\prime}(\tau _{1})|}}g(\tau _{1})e^{itF(\tau _{1})+\frac{i\pi}{4}sgn F^{\prime\prime}(\tau _{1})}
+O(\frac{1}{t}), \hphantom{3a} t\rightarrow \infty .
\end{equation}
In order to compute $F^{\prime\prime}(\tau _{1})$ we use (\ref{sixtwentynine}):
\begin{equation}  \label{sixthirtytwo}
\frac{\partial ^{2}F(\tau _{1},\lambda)}{\partial \tau ^{2}}=\frac{1}{\tau _{1}(1-\tau _{1})}=\frac{(1+\lambda )^{2}}{\lambda }.
\end{equation}
Evaluating equation (\ref{sixten}) at $\tau =\tau _{1}$ we find 
\begin{equation}
F(\tau _{1},\lambda)=\frac{\lambda }{1+\lambda }\ln (\frac{\lambda }{1+\lambda })
+\frac{1 }{1+\lambda }\ln (\frac{1 }{1+\lambda })
+\frac{1 }{1+\lambda }\ln\lambda =-\ln (1+\frac{1}{\lambda }).  \nonumber
\end{equation}
Thus,
\begin{equation}  \label{sixthirtythree}
\sqrt{\frac{2\pi}{t|F^{\prime\prime}(\tau _{1})|}}e^{itF(\tau _{1})+\frac{i\pi}{4}sgn F^{\prime\prime}(\tau _{1})}
=\sqrt{\frac{2\pi}{t}}e^{\frac{i\pi}{4}}\frac{\lambda ^{\frac{1}{2}+it}}{(1+\lambda )^{1+it}}.
\end{equation}
The definition of $G(\sigma ,\tau )$ in (\ref{sixnine}) implies
\begin{equation} \label{sixthirtyfour}
G(\sigma ,\tau _{1})=\bigl( 1-\frac{1}{1+\lambda } \bigr)^{-\frac{1}{2}}\frac{1}{(1+\lambda )^{\sigma -\frac{1}{2}}}=
\frac{\lambda ^{-\frac{1}{2}}}{(1+\lambda )^{\sigma -1}}.
\end{equation}
Furthermore, the definitions of $A(t,\tau )$ and $B(\sigma, t,\tau )$ in (\ref{sixeleven}) and (\ref{sixtwelve}) imply
\begin{equation}  \label{sixthirtyfive}
A(t,\tau _{1})=1+O\left(\frac{1+\lambda }{t\lambda }\right)+O\left(\frac{1+\lambda }{t}\right), \quad t\rightarrow \infty ,
\end{equation}
and
\begin{equation} \label{sixthirtysix}
B(\sigma, t,\tau _{1})=O\left(1,\frac{1}{1+\lambda }, \frac{1}{(1+\lambda )^{2} }, \frac{1}{(1+\lambda )^{3}t^{\sigma }}\right), \quad t\rightarrow \infty .
\end{equation}
Also,
\begin{equation}  \label{sixthirtyseven}
\frac{1}{\tau \pm i\frac{(1-\sigma )}{t}}=\frac{1+\lambda }{1\pm \frac{i(1-\sigma )}{t}(1+\lambda )}.
\end{equation}
Employing equation (\ref{sixthirtyone}) for the estimation of the integrals (\ref{sixthree})-(\ref{sixseven}), and noting that the product of the rhs of equations \eqref{sixthirtythree} and \eqref{sixthirtyfour} equals
\begin{equation*}
\sqrt{\frac{2\pi}{t}} e^{\frac{i\pi}{4}} \frac{\lambda^{it}}{(1+\lambda)^{\sigma+it}},
\end{equation*}
we find equation (\ref{sixtwentyone}), with $P_{j}$, $j=1,\dots,5$ defined in 
(\ref{sixtwentythree})-(\ref{sixtwentyseven}).   

In the above analysis, we have assumed that the stationary points do \textit{not} occur on the boundaries. In particular,
\begin{equation*}
\frac{1}{1+\lambda} \ne 1 - \frac{t^{\delta_3} }{t}, \quad \text{or} \quad \lambda \left(  \frac{t}{t^{\delta_3}} -1 \right) \ne 1.
\end{equation*}
This implies the following constraints:
\begin{align}
\text{for} ~ J^{(1)}&: \quad \frac{m_2}{m_1} \left( \frac{t}{t^{\delta_3}} -1  \right) \ne 1,\\
\text{for} ~ J^{(2)} \text{and} ~ J^{(4)}&: \quad \frac{T}{m} \left( \frac{t}{t^{\delta_3}} -1  \right) \ne 1,\\
\text{for} ~ J^{(3)} \text{and} ~ J^{(5)}&: \quad \frac{m}{T} \left( \frac{t}{t^{\delta_3}} -1  \right) \ne 1.
\end{align}
If a stationary point \textit{does} occur on the boundary, then the relevant contribution is half the contribution computed earlier. We will consider this possible additional term in the analysis presented in lemma \ref{l6.2}.

\textbf{QED}

\vphantom{a}

\begin{corollary}  \label{c6.1}
Let $I_{3}(\sigma ,t,\delta _{2},\delta _{3})$ be defined in (\ref{sixone}). Then,

\begin{equation}
I_{3}(\sigma ,t,\delta _{2},\delta _{3})= \sqrt{\frac{2t}{\pi}}\Re\Biggl(  
e^{-\frac{i\pi}{4}}\sum_{m_{1}=1}^{[T]}\sum_{m_{2}=1}^{[T]}m_{1}^{-\sigma }m_{2}^{-\sigma }
J_{B}^{(1)}(\sigma ,t,\delta _{3},\frac{m_{2}}{m_{1}})  \nonumber
\end{equation}
\begin{equation}
+ \frac{i}{2\pi}  T^{-\sigma }e^{-\frac{i\pi}{4}}\sum_{m=1}^{[T]}m^{-\sigma }
[\frac{i}{2\pi}J_{B}^{(2)}(\sigma ,t,\delta _{3},\frac{T}{m})
+J_{B}^{(4)}(\sigma ,t,\delta _{3},\frac{T}{m})]  \nonumber
\end{equation}
\begin{equation}
+ \frac{i}{2\pi} T^{-\sigma }e^{-\frac{i\pi}{4}}\sum_{m=1}^{[T]}m^{-\sigma }
[-\frac{i}{2\pi}J_{B}^{(3)}(\sigma ,t,\delta _{3},\frac{m}{T})
+J_{B}^{(5)}(\sigma ,t,\delta _{3},\frac{m}{T})]  \nonumber
\end{equation}
\begin{equation}
+T^{-2\sigma }e^{-\frac{i\pi}{4}}
J^{(6)}(\sigma ,t,\delta _{2},\delta _{3}) \Biggr) \left[1+O\left( \frac{1}{t}\right) \right] 
 \nonumber
\end{equation}
\begin{equation}
+ 2 \Re\Biggl(
\sum_{m_{1}=1}^{[T]}\sum_{m_{2}=1}^{[T]}\frac{1}{(m_{1}+m_{2})^{s}m_{2}^{\bar{s} }}
\left[1 + O\left(\frac{1+m_1/m_2}{t}\right) + O\left(\frac{1+m_2/m_1}{t}\right) \right] 
 \nonumber
\end{equation}
\begin{equation}
+ \frac{i}{2\pi}  \sum_{m=1}^{[T]}\frac{1}{T^{\bar{s}}(m+T)^{s} }
\Biggl( 
\frac{1+\frac{T}{m}}{1-\frac{i(1-\sigma )}{t}(1+\frac{T}{m})} 
\nonumber
\end{equation}
\begin{equation}
+O\bigl(  1,\frac{1}{1+\frac{T}{m}},\frac{1}{(1+\frac{T}{m})^{2}},\frac{1}{(1+\frac{T}{m})^{3}t^{\sigma }}  \bigr)
\left( 1+O\left(\frac{1+T/m}{t}\right) +O\left(\frac{1+m/T}{t}\right) \right)
\Biggr)  \nonumber
\end{equation}
\begin{equation}
+ \frac{i}{2\pi}  \sum_{m=1}^{[T]}\frac{1}{m^{\bar{s}}(m+T)^{s} }
\Biggl( 
\frac{1+\frac{m}{T}}{1+\frac{i(1-\sigma )}{t}(1+\frac{m}{T})} 
\nonumber
\end{equation}
\begin{equation}
+O\bigl(  1,\frac{1}{1+\frac{m}{T}},\frac{1}{(1+\frac{m}{T})^{2}},\frac{1}{(1+\frac{m}{T})^{3}t^{\sigma }}  \bigr)\left( 1+O\left(\frac{1+T/m}{t}\right) +O\left(\frac{1+m/T}{t}\right) \right)
\Biggr)  \nonumber
\end{equation}
\begin{equation}  \label{sixthirtyeight}
\times \left[1+O\left(\frac{1}{t}\right)\right],  \quad  \frac{1}{2}\le\sigma <1, \quad t\rightarrow \infty ,
\end{equation}
where the integrals $\{ J_{B}^{(j)} \}_{1}^{5}$ are obtained from the integrals $\{ J^{(j)} \}_{1}^{5}$ defined in 
(\ref{sixthree})-(\ref{sixseven}) by replacing the contour of integration with the ray from the point $1-t^{\delta _{3}-1}$ to 
$\infty e^{i\varphi }$ with $\varphi $ satisfying (\ref{sixtwentytwo}) and the integral 
$J^{(6)}(\sigma ,t,\delta _{2},\delta _{3} )$ is defined in (\ref{sixeight}). The integers $(m_{1},m_{2})$ 
in the double sum, as well as the integer $m$ in the first and the second single sum, satisfy the inequalities
\begin{equation} \label{sixthirtynine}
\frac{1}{t^{1-\delta _{3}}-1}<\frac{m_{2}}{m_{1}}<t^{1-\delta _{2}}-1, \hphantom{3a} 
m>\frac{T}{t^{1-\delta _{2}}-1}, \hphantom{3a} 
m>\frac{T}{t^{1-\delta _{3}}-1},
\end{equation}
respectively.
\end{corollary}

\textbf{Proof} We replace in (\ref{sixtwo}) the expressions of $J^{(1)}$, $J^{(2)}$, $J^{(3)}$, $J^{(4)}$, $J^{(5)}$ with the expressions of the rhs of (\ref{sixtwentyone}) with $\lambda $ given for $J^{(1)}$, $\left\{ J^{(2)}, J^{(4)}\right\}$,  $\left\{ J^{(3)}, J^{(5)}\right\}$, by $m_{2}/m_{1}$, $T/m$, $m/T$ respectively. In this connection the following formulae are valid \\

\noindent $J^{(1)}$:
\begin{equation}
m_{1}^{-\sigma }m_{2}^{-\sigma }\frac{\lambda ^{it}}{(1+\lambda )^{\sigma +it}}
\left[ 1+O\left(\frac{1+\frac{1}{\lambda}}{t}\right) +O\left(\frac{1+\lambda}{t} \right)\right]\Big|_{\lambda =\frac{m_{2}}{m_{1}}} \nonumber
\end{equation}
\begin{equation}  \label{sixfourty}
=\frac{1}{m_{2}^{\bar{s}}(m_{1}+m_{2})^{s} }
\left[1+O\left(\frac{1}{t}\left(1+\frac{m_{1}}{m_{2}}\right)\right)+O\left(\frac{1}{t}\left(  1 + \frac{m_2}{m_1} \right) \right) \right].
\end{equation}
$J^{(2)}$:
\begin{multline}  \label{sixfourtyone}
(\frac{i}{2\pi})T^{-\sigma }m^{-\sigma }\frac{\lambda ^{it}}{(1+\lambda )^{\sigma +it}}
\left[ 1+O\left(\frac{1+\frac{1}{\lambda}}{t}\right) +O\left(\frac{1+\lambda}{t} \right)\right] \\
\times\frac{(1+\lambda )}{1-\frac{i(1-\sigma )}{t}(1+\lambda )}\Big|_{\lambda =\frac{T}{m}} \\ 
= \frac{i}{2\pi}\frac{1}{T^{\bar{s}}m(m+T)^{-1+s} }
\frac{\left[1+O\left(\frac{1+m/T}{t}\right)+O\left(\frac{1+T/m}{t}\right)\right]}{1- \frac{i(1-\sigma )}{t}(1+\frac{T}{m})}.
\end{multline}
$J^{(3)}$:
\begin{multline}  \label{sixfourtytwo}
(-\frac{i}{2\pi})T^{-\sigma }m^{-\sigma }\frac{\lambda ^{it}}{(1+\lambda )^{\sigma +it}}
\left[ 1+O\left(\frac{1+\frac{1}{\lambda}}{t}\right) +O\left(\frac{1+\lambda}{t} \right)\right] \\
\times \frac{(1+\lambda )}{1+\frac{i(1-\sigma )}{t}(1+\lambda )}
\mid _{\lambda =\frac{m}{T}}  \\
= - \frac{i}{2\pi}\frac{1}{Tm^{\bar{s}}(m+T)^{-1+s} }
\frac{\left[1+O\left(\frac{1+m/T}{t}\right)+O\left(\frac{1+T/m}{t}\right)\right]}{1+\frac{i(1-\sigma )}{t}(1+\frac{m}{T})}.
\end{multline}
$J^{(4)}$:
\begin{multline} \label{sixfourtythree}
T^{-\sigma }m^{-\sigma }\frac{\lambda ^{it}}{(1+\lambda )^{\sigma +it}}
\left[ 1+O\left(\frac{1+\frac{1}{\lambda}}{t}\right) +O\left(\frac{1+\lambda}{t} \right)\right] \\ 
\times O\left( 1,\frac{1}{1+\lambda }, \frac{1}{(1+\lambda )^{2} }, \frac{1}{(1+\lambda )^{3}t^{\sigma }}  \right)
\mid _{\lambda =\frac{T}{m}}  \\
=\frac{1}{T^{\bar{s}}(m+T)^{s} } \left[ 1+O\left(\frac{1+m/T}{t}\right)+O\left(\frac{1+T/m}{t}\right) \right] \\
\times O\left( 1,\frac{1}{1+\frac{T}{m} }, \frac{1}{(1+\frac{T}{m} )^{2} }, \frac{1}{(1+\frac{T}{m} )^{3}t^{\sigma }}  \right).
\end{multline}
$J^{(5)}$:
\begin{multline} \label{sixfourtyfour}
T^{-\sigma }m^{-\sigma }\frac{\lambda ^{it}}{(1+\lambda )^{\sigma +it}}
\left[ 1+O\left(\frac{1+\frac{1}{\lambda}}{t}\right) +O\left(\frac{1+\lambda}{t} \right)\right] \\
\times O\left( 1,\frac{1}{1+\lambda }, \frac{1}{(1+\lambda )^{2} }, \frac{1}{(1+\lambda )^{3}t^{\sigma }}  \right)
\mid _{\lambda =\frac{m}{T}} \\
=\frac{1}{m^{\bar{s}}(m+T)^{s} } \left[ 1+O\left(\frac{1+m/T}{t}\right)+O\left(\frac{1+T/m}{t}\right) \right] \\
\times O\left( 1,\frac{1}{1+\frac{m}{T} }, \frac{1}{(1+\frac{m}{T} )^{2} }, \frac{1}{(1+\frac{m}{T} )^{3}t^{\sigma }}  \right) .
\end{multline}
The inequalities (\ref{sixtwentyeight}) with $\lambda $ given by
\begin{equation}
\lambda =\frac{m_{2}}{m_{1}}, \quad \lambda =\frac{T}{m}, \quad \lambda =\frac{m}{T},  \nonumber
\end{equation}
yield the inequalities (\ref{sixthirtynine}) respectively.

\textbf{QED}

In what follows we analyse the contributions of the integrals $\{ J_{B}^{(j)}(\sigma, t, \delta, \lambda)\}_1^5$ and $J^{(6)}$.

\begin{lemma} \label{l6.2}
Define the integrals $\{ J_{B}^{(j)}(\sigma, t, \delta, \lambda)\}_1^5$ by
\begin{multline} \label{sixLT1}
 J_{B}^{(j)}(\sigma, t, \delta, \lambda) = \int_{1-t^{\delta-1}}^{\infty e^{i\varphi}}  (1-z)^{-\frac{1}{2}}z^{\sigma-\frac{1}{2}}e^{itF(z,\lambda)} D^{(j)}(\sigma, t, z) \textrm{d}z, \\
j=1,2,3,4,5, \quad \frac{1}{2}\le\sigma<1, \quad t>0,
\end{multline}
where $\delta$ is a sufficiently small, positive constant, $\varphi$ is a positive constant satisfying the inequality \eqref{sixtwentytwo}, the function $F$ is defined by \eqref{sixten}, the functions $D^{(j)}$ are given for $j=1,2,3,4$ by
\begin{equation} \label{sixLT2}
A(t,z), \quad \frac{A(t,z)}{z - \frac{i(1-\sigma)}{t}}, \quad \frac{A(t,z)}{z + \frac{i(1-\sigma)}{t}}, \quad A(t,z) \overline{B(\sigma,t,z)}, \quad A(t,z)B(\sigma,t,z),
\end{equation}
respectively with $A$ and $B$ defined in \eqref{sixeleven} and \eqref{sixtwelve}, and $\lambda$ given for $j=1$ by $m_2/m_1$, for $j=2$ and 4 by $T/m$, and for $j=3$ and 5 by  $m/T$.

Define
$J^{(6)}(\sigma,t,\delta)$ by the rhs of \eqref{sixLT1} but with $F(\tau,\lambda)$ replaced with $F(\tau,1)$. Let $\lambda=\frac{m_2}{m_1}$.

If $(m_1,m_2)\in N$, then
\begin{align} \label{sixLT3}
 J_{B}^{(j)}(\sigma, t,\delta, \lambda) = 
\frac{i\lambda^{i(t-t^{\delta})}}{t^{\frac{\delta+1}{2}}}
& D^{(j)}\left(\sigma,t,1-t^{\delta-1}\right)
\dfrac{t^{i(\delta-1)t^{\delta}}(1-t^{\delta-1})^{\sigma-\frac{1}{2}
+i(t-t^{\delta})}}{\ln \left[ \lambda  \left( t^{1-\delta } - 1 \right) \right]} \nonumber\\
&+ o\left( \dfrac{1}{t^{\frac{\delta+1}{2}} \ln \left[ \lambda  \left( t^{1-\delta } - 1 \right) \right] }\right), \qquad t\rightarrow\infty.
\end{align}
If $(m_1,m_2)\in N^c$, then
\begin{align} \label{sixLT4}
 J_{B}^{(j)}(\sigma, t,\delta, \lambda) = 
\sqrt{\frac{2}{t}} \lambda^{i(t-t^{\delta})} &
D^{(j)}\left(\sigma,t,1-t^{\delta-1}\right)\big(\Omega +o(1)\big) \nonumber \\ &+  o\left(t^{-\frac{1}{2}}\right), \qquad t\rightarrow\infty,
\end{align}
where $\Omega$ is complex finite constant and  the set $N$ is defined in (1.13b).
\end{lemma}
\textbf{Proof.}

%

Let $F(z,\lambda)$ be defined by \eqref{sixten}. Then, the identity
$$\dfrac{\partial F}{\partial z}(z,\lambda)=\ln\left(\frac{z\lambda}{1-z}\right)$$
implies that the integrals $\left\{ J_{B}^{(j)} \right\}_1^5$ have a stationary point at $z=\frac{1}{1+\lambda}$. If 
\begin{equation}\label{sixLTL1}
\frac{1}{t^{1-\delta}-1}< \lambda \leq t^{\delta-1}-1,
\end{equation} 
the stationary point is in the interval $\left( 0,1-t^{\delta-1}\right)$, and thus it is away from the contour of integration. However, when $\lambda=\lambda_c$, where 
\begin{equation}\label{sixLTL2}
\lambda_{c}=\frac{t^{\delta-1}}{1-t^{\delta-1}},
\end{equation}
then the stationary point is at $z=1-t^{\delta-1}$, i.e. at the left endpoint of integration.

In order to analyse the behavior of $J_{B}$ for
$\lambda$ near $\lambda_{c}$, we introduce new variables
$\Lambda$ and $\zeta$ such that $\Lambda=0$ and $\zeta=0$
correspond to $\lambda=\lambda_{c}$ and $z=1-t^{\delta-1}$:
\begin{equation}\label{sixLTLa}
\lambda=\lambda_{c}(1+\Lambda),~~~~z=\frac{1}{1+\lambda_{c}}(1+\lambda_{c}\zeta).
\end{equation}
Then, $1-z=\lambda_{c}(1-\zeta)/(1+\lambda_{c})$, thus,
\begin{equation}\label{sixLT5}
J_{B}^{(j)}=\left(\frac{\lambda_{c}}{1+\lambda_{c}}\right)^{\frac{1}{2}}
\left(\frac{1}{1+\lambda_{c}}\right)^{\sigma-\frac{1}{2}}
e^{\frac{itf_0(\Lambda,\lambda_c)}{1+\lambda_{c}}}\tilde{J}_{B}^{(j)}, \quad j=1,\ldots,5,
\end{equation}
where $\tilde{J}_{B}^{(j)}$ is given by
\begin{align}\label{sixLT6}
\tilde{J}_{B}^{(j)}=\int_{0}^{\infty e^{i\varphi}}(1-\zeta)^{-\frac{1}{2}}
(1+\lambda_{c}\zeta)^{\sigma-\frac{1}{2}}
e^{\frac{itf_1(\zeta,\Lambda,\lambda_c)}{1+\lambda_{c}}} & D^{(j)}\left(\sigma,t,\frac{1+\lambda_{c}\zeta}{1+\lambda_{c}}\right) d\zeta, \nonumber\\
&j=1,\ldots,5,
\end{align}
with
\begin{equation}\label{sixLTdf0}
f_0=\ln\left[\frac{\lambda_{c}(1+\Lambda)}{1+\lambda_{c}}\right]
-\lambda_{c}\ln \left(\frac{1+\lambda_{c}}{\lambda_{c}}\right),
\end{equation}
and
\begin{equation}\label{sixLTdf1}
f_1=\lambda_{c}\zeta[\ln(1+\lambda_{c}\zeta)-\ln (1-\zeta)+\ln(1+\Lambda)]+\ln(1+\lambda_{c}\zeta)+\lambda_{c}\ln (1-\zeta).
\end{equation}

Rigorous uniform asymptotics of $\left\{ \tilde{J}_{B}^{(j)} \right\}_1^6$, as $t \to \infty$ is derived in \cite{FFS}:

\textbf{Theorem \cite{FFS}} Define $\omega$ by
\begin{equation}\label{sixLTom}
\omega(t,\Lambda)=\sqrt{\frac{\lambda_{c}t}{2}}\frac{\ln (1+\Lambda)}{1+\lambda_{c}}.
\end{equation}
Observe that $\omega\geq 0$ since $\Lambda\geq 0$.
The leading-order asymptotics of $J_B$ is given by \eqref{sixLT5}, with the leading-order asymptotics of $\tilde{J}_B$ given by
\begin{equation}\label{eq:result1}
\tilde{J}_B=\sqrt{\frac{2}{\lambda_c t}}e^{-i\omega^{2}} \left( \int_\omega^{\infty e^{i \pi/4}}e^{i \xi^2}d \xi\right) \left(1+ o(1)\right),
\end{equation}
where the $o(1)$ term is independent of $\Lambda$.

\underline{\textbf{QED}}

\begin{remark}
The above theorem is proven in \cite{FFS} in the particular case of $D^{(1)}$. However, it is straightforward to extend the proof to the case that  $D^{(1)}$ is replaced by $\left\{ D^{(j)}\right\}_2^5$.
\end{remark}

\begin{remark}
If $\Lambda$ is such that $\omega=\mathcal {O}(1)$ as $t\rightarrow\infty$ (for example~$\Lambda=0$), then the integral on the right-hand side of \eqref{eq:result1}  is an $\mathcal {O}(1)$ quantity independent of $\Lambda$.
If $\Lambda$ is such that $\omega\rightarrow\infty$ as $t\rightarrow\infty$, then
\begin{equation}\label{eq:result2}
\tilde{J}_B=
\sqrt{\frac{2}{\lambda_c t}}
\left( \frac{-1}{2i \omega}+ \mathcal {O}\left( \frac{1}{\omega^3}\right)\right) \bigg(1+ o(1)\bigg),~~t\rightarrow\infty,
\end{equation}
 where both the $o(1)$ and the omitted constant in the $\mathcal {O}( 1/\omega^3)$ term
are independent of $\Lambda$.
\end{remark}

{\remark The set $N$ is the complement of the set defined, for $\lambda=\frac{m_2}{m_1}$, via the condition that $\omega=O(1)$.  Thus, letting $\mathcal{P} =\frac{m_2}{m_1} \left(t^{1-\delta} -1\right)$ and using equations  \eqref{sixLTL2}, \eqref{sixLTLa} and \eqref{sixLTom}, we  obtain that $\frac{1}{\ln \mathcal{P}}\ll t^{\delta/2}$. Also, the restrictions $m_2\leq T$ and $m_1\geq 1$ imply $\frac{m_2}{m_1}< t^2,$ thus $\frac{1}{\ln \mathcal{P}}> \frac{1}{2\ln t}.$
 }

\begin{corollary} \label{c6.2}
Let $I_{3}(\sigma ,t,\delta _{2},\delta _{3})$ be defined in (\ref{sixone}). Then,

\begin{align*}
I_{3}(\sigma ,t,\delta _{2},\delta _{3})&=  \frac{1}{\pi} \sqrt{\frac{2}{t^{\delta}}}  \Re \vast\{ e^{\frac{i\pi}{4}} \left( t^{\delta_3-1}\right)^{it^{\delta_3}} \left(  1 -t^{\delta_3-1} \right)^{\sigma- \frac{1}{2}+i(t-t^{\delta_3})}   \\
& \left.  \times  \mathop{\sum\sum}_{m_1,m_2 \in N}\frac{1}{m_1^{s-it^{\delta_3}}}  \frac{1}{m_2^{\bar{s} + it^{\delta_3}}} \frac{1}{\ln{\left[ \frac{m_2}{m_1} \left(  t^{1-\delta_3} - 1 \right)  \right]}}   \left[ 1 + O\left( t^{-\delta_3}\right) + o(1)  \right] \right\} 
\end{align*}
\begin{align*}
+ \frac{1}{\sqrt{2}\pi^{3/2}} \frac{1}{t^\frac{\delta_3}{2}}\frac{1}{T^\sigma}  &\Re \vast\{ e^{i\frac{3\pi}{4}} \left( t^{\delta_3-1}\right)^{it^{\delta_3}} \left(  1 -t^{\delta_3-1} \right)^{\sigma- \frac{1}{2}+i(t-t^{\delta_3})}  \\
& \left.  \times  \sum_{m=1}^{[T]}\frac{1}{m^{s-it^{\delta_3}}}   \frac{1}{\ln{\left[ \frac{T}{m} \left(  t^{1-\delta_3} - 1 \right)  \right]}} [1+o(1)] \right\} \times   O\left(1, 1-t^{\delta_3-1}, t^{-\delta_3}\right) 
\end{align*}
\begin{align*}
+ \frac{1}{\sqrt{2}\pi^{3/2}} \frac{1}{t^\frac{\delta_3}{2}}\frac{1}{T^\sigma}  &\Re \vast\{ e^{i\frac{3\pi}{4}} \left( t^{\delta_3-1}\right)^{it^{\delta_3}} \left(  1 -t^{\delta_3-1} \right)^{\sigma- \frac{1}{2}+i(t-t^{\delta_3})}  \\
& \left.  \times  \sum_{m\in \tilde{N}}\frac{1}{m^{s-it^{\delta_3}}}   \frac{1}{\ln{\left[ \frac{m}{T} \left(  t^{1-\delta_3} - 1 \right)  \right]}} [1+o(1)] \right\} \times   O\left(1, 1-t^{\delta_3-1}, t^{-\delta_3}\right) 
\end{align*}
\begin{equation}
+ 2 \Re\Bigg\{
\mathop{\sum\sum}_{m_1,m_2 \in M}\frac{1}{(m_{1}+m_{2})^{s}m_{2}^{\bar{s} }}
\left[1 + O\left(\frac{1+m_1/m_2}{t}\right) + O\left(\frac{1+m_2/m_1}{t}\right) \right] 
 \nonumber
\end{equation}
\begin{equation}
+ \frac{i}{2\pi}  \sum_{m \in \tilde{M}}\frac{1}{T^{\bar{s}}(m+T)^{s} }
\Biggl( 
\frac{1+\frac{T}{m}}{1-\frac{i(1-\sigma )}{t}(1+\frac{T}{m})} 
\nonumber
\end{equation}
\begin{equation}
+O\bigl(  1,\frac{1}{1+\frac{T}{m}},\frac{1}{(1+\frac{T}{m})^{2}},\frac{1}{(1+\frac{T}{m})^{3}t^{\sigma }}  \bigr)
\left( 1+O\left(\frac{1+T/m}{t}\right) +O\left(\frac{1+m/T}{t}\right) \right)
\Biggr)  \nonumber
\end{equation}
\begin{equation}
+ \frac{i}{2\pi}  \sum_{m \in \hat{M}}\frac{1}{m^{\bar{s}}(m+T)^{s} }
\Biggl( 
\frac{1+\frac{m}{T}}{1+\frac{i(1-\sigma )}{t}(1+\frac{m}{T})} 
\nonumber
\end{equation}
\begin{equation}
+O\bigl(  1,\frac{1}{1+\frac{m}{T}},\frac{1}{(1+\frac{m}{T})^{2}},\frac{1}{(1+\frac{m}{T})^{3}t^{\sigma }}  \bigr)\left( 1+O\left(\frac{1+T/m}{t}\right) +O\left(\frac{1+m/T}{t}\right) \right)
\Bigg\}  \nonumber
\end{equation}
\begin{equation}  \label{sixthirtyeight-new}
\times \left[1+O\left(\frac{1}{t}\right)\right] \ + \begin{cases}O\left(t^\frac{\delta_3}{2} \ln t\right), \quad &\sigma=\frac{1}{2}, \\ O\left(t^{\frac{1}{2}-\sigma+\delta_3\sigma}\right),   &\frac{1}{2}<\sigma <1, \end{cases} \quad t\rightarrow \infty .
\end{equation}
The set $N$ is defined in \eqref{1dN} and  
\begin{multline}  \label{1dNt}
\tilde{N}= \left\{ m\in\mathbb{N}^{+}, \ 1\leq m\leq [T], \  m>\frac{T}{t^{1-\delta _{3}}-1}\bigg(1+c(t)\bigg), \ t^{-\frac{\delta_3}{2}}\ll c(t) \ll 1\right\}
\end{multline}

Furthermore, the set $M$ is defined in \eqref{1.8},
\begin{equation} \label{1dMt}
\tilde{M}= \left\{ m\in\mathbb{N}^{+}, \ 1\leq m\leq [T], \  m>\frac{T}{t^{1-\delta _{2}}-1}\right\}
\end{equation}
and
\begin{equation}\label{1dMh}
\hat{M}= \left\{ m\in\mathbb{N}^{+}, \ 1\leq m\leq [T], \  m>\frac{T}{t^{1-\delta _{3}}-1}\right\}.
\end{equation}
\end{corollary}
\textbf{Proof} We employ the results of Lemma \ref{l6.2} in the terms $ J_B^{(j)}(\sigma,t,\delta,\lambda), \ j=1,\ldots,5,$ of Corollary \ref{c6.1}. If $\lambda\in N$, we obtain the first three terms and the proof is given by following the steps of the proof of Corollary \ref{c6.1}, but now the terms $D^{(j)}$ are much simpler than the respective terms of Corollary \ref{c6.1}. The set $\tilde{N}$ is the same with the set $N$, with the substitution $\frac{m_2}{m_1}\to \frac{m}{T}$.

It is straightforward to show using the techniques employed in section \ref{sec7}, that if $\lambda \in N^c$, then the contribution of $ J_B^{(j)}, \ j=1,\ldots,5,$ yields a  sum which is of order $O\left(t^\frac{\delta_3}{2} \ln t \right)$ if $\sigma=\frac{1}{2}$, and of order $O\left(t^{\frac{1}{2}-\sigma+\delta_3\sigma}\right)$ if $\frac{1}{2}<\sigma<1$, as $t\rightarrow\infty$. In particular, the leading order term of this sum is the same with the sum $S_A$ in Lemma \ref{lemma7.4}. Thus, the relevant estimate can be obtained in the same way as the estimate derived in Lemma \ref{lemma7.4}.

Finally, the contribution from $J^{(6)}$ is of order $O(t^{-2\sigma -\frac{1}{2} -\frac{\delta}{2} }( \ln{t})^{-1})$, thus neglected. 

\underline{\textbf{QED}}\\

\begin{remark} The $o(1)$  term appearing in Lemma \ref{l6.2} and in Corollary \ref{c6.2} contains terms which are functions of $\lambda$, namely these terms depend on $\frac{m_2}{m_1}, \frac{T}{m}, \frac{m}{T}$, and they are much smaller than the leading order term.
\end{remark}

\begin{remark} A more detailed analysis for the contribution of $ J_B^{(j)}, \ j=1,\ldots,5,$ when $\lambda \in N^c$, is presented in \cite{FK}, where the relevant sum is treated using techniques that appear in \cite{T} and \cite{T2}. In particular, the oscillatory part of the sum is of the form $e^{i\left[t\ln\lambda + O(t^\delta )\right]}$, and hence these techniques can be indeed applied.

At the limiting case where the stationary points are on the end-point of the interval of integration, this analysis simplifies; see Remark \ref{rem6.1} for the details.
\end{remark}

\begin{theorem} \label{t6.1}
Define the set $M$ by \eqref{1.8}, i.e.,
\begin{multline}  \label{sixfortysix}
M = \left\{ m_{1}\in\mathbb{N}^{+}, \quad  m_{2}\in\mathbb{N}^{+}, \quad 1\leq m_{1}\leq [T], \quad 1\leq m_{2}< [T], \right.  \\
\left. \frac{1}{t^{1-\delta _{3}}-1}<\frac{m_{2}}{m_{1}}<t^{1-\delta _{2}}-1, \quad  t>0, \quad  T=\frac{t}{2\pi} \right\}
\end{multline}
and the set $N$ by \eqref{1dN}, i.e.,
\begin{multline*}
N= \bigg\{ m_{1}\in\mathbb{N}^{+}, \quad  m_{2}\in\mathbb{N}^{+}, \quad 1\leq m_{1}\leq [T], \quad 1\leq m_{2}< [T],  \\
 \frac{m_{2}}{m_{1}}>\frac{1}{t^{1-\delta _{3}}-1}\bigg(1+c(t)\bigg), \ t^{-\frac{\delta_3}{2}}\ll c(t) \ll 1, \quad  t>0, \quad  T=\frac{t}{2\pi} \bigg\},
\end{multline*}
where $\delta _{2}$ and $\delta _{3}$ are sufficiently small, positive constants.

Let $I_3(\sigma,t,\delta_2,\delta_3)$, $\ \frac{1}{2}\le\sigma<1$, $t>0$, be defined in \eqref{sixone}. Then,
\begin{align}  \label{sixfortyseven}
 I_{3}&(\sigma ,t,\delta _{2},\delta _{3})=    \frac{1}{\pi} \sqrt{\frac{2}{t^{\delta_3}}}  \Re \vast\{ e^{\frac{i\pi}{4}} \left( t^{\delta_3-1}\right)^{it^{\delta_3}} \left(  1 -t^{\delta_3-1} \right)^{\sigma- \frac{1}{2}+i(t-t^{\delta_3})} \nonumber \\
&\times  \mathop{\sum\sum}_{m_1,m_2 \in N}\frac{1}{m_1^{s-it^{\delta_3}}}  \frac{1}{m_2^{\bar{s} + it^{\delta_3}}} \frac{1}{\ln{\left[ \frac{m_2}{m_1} \left(  t^{1-\delta_3} - 1 \right)  \right]}}   \left[ 1 + O\left( t^{-\delta_3}\right) + o(1)  \right] \vast\} \nonumber  \\
+2&  \Re{ \left\{  \mathop{\sum\sum}_{m_1,m_2 \in M} \frac{1}{(m_1+m_2)^s m_2^{\bar{s}}}  \right\}}+ \begin{cases}O\left(t^\frac{\delta_3}{2} \ln t\right), \quad &\sigma=\frac{1}{2}, \\ O\left(t^{\frac{1}{2}-\sigma+\delta_3\sigma}\right),   &\frac{1}{2}<\sigma <1, \end{cases}\nonumber  \\
& \qquad + O \left( t^{1-2\sigma} \ln{t} \right), \qquad t \to \infty.
\end{align}
\end{theorem}

\textbf{Proof} Equation \eqref{sixfortyseven} is a direct consequence of equation \eqref{sixthirtyeight-new} and of the estimation of the error terms associated with the contribution of the stationary points and the end point. Regarding the former terms, we note that the first sum of the stationary points contribution occurring in \eqref{sixthirtyeight-new} is given by
\begin{equation*}
\sum_{m\in\tilde{M}}\frac{1+\frac{T}{m}}{T^{\bar{s}}(m+T)^{s} } \left[ 1+O \left(\frac{1}{m}, \frac{1}{1+\frac{T}{m}} \right) \right].
\end{equation*}
Since $m\in\tilde{M}$, we find that
\begin{equation*}
m\geq\frac{T}{t^{1-\delta _{2}}-1} > \frac{t/2\pi}{t^{1-\delta _{2}}} = \frac{t^{\delta_2}}{2\pi},
\end{equation*}
and then we can employ the following ``crude'' estimates:
\begin{multline} \label{newsixfortyseven}
\left| \sum_{\frac{t^{\delta _{2}}}{2\pi}}^{[T]}\frac{1}{T^{\bar{s}}(m+T)^{s} }\left(1+\frac{T}{m}\right)\right|  \leq
\int_{\frac{t^{\delta _{2}}}{2\pi}}^{T}\frac{\left(1+\frac{T}{x}\right)}{T^{\sigma }(x+T)^{\sigma }}{\text{d}}x \\
=T^{1-2\sigma }
\int_{\frac{t^{\delta _{2}}}{2\pi}}^{T}\frac{\left(1+\frac{T}{x}\right)^{1-\sigma }}{(\frac{x}{T})^{\sigma }}{\text{d}}\left(\frac{x}{T}\right)=
T^{1-2\sigma }
\int_{t^{\delta _{2}-1}}^{1}\frac{(1+\rho)^{1-\sigma }}{\rho}{\text{d}}\rho =
O(T^{1-2\sigma } \ln{t}). 
\end{multline}
Indeed, for $t^{\delta _{2}-1}<\rho<1$, $(1+\rho)^{1-\sigma }$ is bounded. Thus,
\begin{equation*}
O\left(  \int_{t^{\delta _{2}-1}}^{1}\frac{(1+\rho)^{1-\sigma }}{\rho}d\rho \right) = O\left (  \int_{t^{\delta _{2}-1}}^{1}\frac{d\rho}{\rho}  \right) = O(\ln{t}).
\end{equation*}
Thus, the highest order term of this sum is of order $O(t^{1-2\sigma}\ln{t})$.

For the associated first error term which involves $1/m$, we can use partial summation, the fact that $m>t^{\delta_2}/2\pi$, and \eqref{newsixfortyseven}, to show that this term is of order $O( t^{1 - 2\sigma - \delta_2}\ln{t} )$.

For the associated second error term, we have to estimate the sum
\begin{equation*}
\sum_{m\in\tilde{M}}\frac{1}{T^{\bar{s}}(m+T)^{s} }.
\end{equation*}
Employing \eqref{newsixfortyseven} we find that this term is of order $O(t^{1-2\sigma})$.

The second single sum of the stationary points contribution occurring in \eqref{sixthirtyeight-new} is given by
\begin{equation*}
\sum_{m\in\hat{M}}\frac{1+\frac{m}{T}}{m^{\bar{s}}(m+T)^{s} } \left[ 1+O \left(\frac{1}{m}, \frac{1}{1+\frac{m}{T}} \right) \right].
\end{equation*}
In analogy with \eqref{newsixfortyseven}, it is straightforward to show that
\begin{equation}
\left| \sum_{\frac{t^{\delta _{3}}}{2\pi}}^{[T]}\frac{1}{m^{\bar{s}}(m+T)^{s} } \left( 1+\frac{m}{T} \right) \right| \leq \int_{\frac{t^{\delta _{3}}}{2\pi}}^{T}\frac{ \left( 1+\frac{x}{T} \right)}{x^{\sigma }(x+T)^{\sigma }}{\text{d}}x= 
O(T^{1-2\sigma }).
\end{equation}
Proceeding as earlier, we find that the highest order term is of order $O( t^{1-2\sigma})$, the first error term is of order $O( t^{1-2\sigma -\delta_3} )$, and the second error term is of order $O( t^{1-2\sigma})$.


It is worth noting that the above single sums involve the term $\exp { \left\{ if(T,m) \right\} }$, where
\begin{subequations} \label{sixfortyfive}
\begin{equation} \label{sixfortyfivea}
f(T,m)= 2 \pi T \ln{\left( 1 +\frac{T}{m}\right)}, \quad t>0, \quad m\in\mathbb{Z}^{+},
\end{equation}
or
\begin{equation} \label{sixfortyfiveb}
f(T,m)= 2 \pi T \ln{\left( 1 +\frac{m}{T}\right)}, \quad t>0, \quad m\in\mathbb{Z}^{+}.
\end{equation}
\end{subequations}
If $f(T,m)$ is given by \eqref{sixfortyfivea}, it is straightforward to show that the basic properties of $f(T,m)$ are similar with the properties of the analogous exponential occurring in the Riemann zeta function. If $f(T,m)$ is given by \eqref{sixfortyfiveb}, the situation is slightly more complicated. In both cases, it is possible to apply the classical techniques to estimate the relevant sums:\\ \\
(a) For the sum involving $m^{-\bar{s}}(m+T)^{-s}(1+\frac{m}{T})$, which corresponds to \eqref{sixfortyfivea}, the techniques of \cite{T} yield an estimate which is better than the ``crude'' estimate obtained above. Indeed, for the highest order term, as well as for the second error term we find the following estimates:
\begin{align*}
\sigma = \frac{1}{2}&: \quad O \left(  t^{-\frac{1}{3}} \ln{t} \right), \\
\sigma = 1&: \quad O \left(  t^{-1} \ln{t} \right).
\end{align*}
Thus, the Phragm\'en-Lindel{\"o}f convexity principle implies that for $\frac{1}{2}\le\sigma<1$, the relevant term is of order $O(t^{\frac{1}{3} - \frac{4}{3} \sigma }\ln{t})$. Similarly, for the first error term we find the following estimates:
\begin{align*}
\sigma = \frac{1}{2}&: \quad O \left(  t^{-\frac{1}{2}} \ln{t} \right), \\
\sigma = 1&: \quad O \left(  t^{-1} \right).
\end{align*}
Thus, the relevant term is of order $O(t^{-\sigma})$, $\frac{1}{2}\le\sigma<1$.
\\ \\
(b) For the sum involving $T^{-\bar{s}}(m+T)^{-s}(1+\frac{T}{m})$, which corresponds to \eqref{sixfortyfiveb}, the techniques of \cite{T} give the estimate $O (t^{1-2\sigma}\ln{t})$, which is exactly the estimate obtained via the ``crude'' estimates. In this case, we only find an improvement for the second error term of the first sum, where we obtain $O(t^{ \frac{1}{2} - 2\sigma }\ln{t})$ instead of $O(t^{ 1 - 2\sigma})$, see \cite{FK} for details.

The first sum of the end point contribution, i.e. the first single sum of \eqref{sixthirtyeight-new}, gives the following contribution:
\begin{equation*}
\frac{1}{t^{\delta_3 /2}} \frac{1}{T^{\sigma}} \sum_{m=1}^{[T]} m^{-\sigma} m^{-i (t-t^{\delta_3 })} \frac{1}{ \ln{ \left[  \frac{T}{m} (t^{1-\delta_3 }-1) \right]  }}   \times O \left( 1, 1 - t^{\delta_3  - 1}, t^{-\delta_3 } \right).
\end{equation*}
Since $m\in(1,T)$, it follows that there exist positive constants $A_1$ and $A_2$ such that
\begin{equation*}
A_1\ln{t} < \ln{  \left[ \frac{T}{m} (t^{1-\delta_3  } - 1)  \right] } < A_2\ln{t}.
\end{equation*}
Then, the estimate \eqref{newsixfortyseven} implies that this term is of order $O(t^{-\frac{\delta_3}{ 2}}T^{1-2\sigma}(\ln{t})^{-1})$.

It is possible to improve this estimate by using the classical techniques together with partial summation in order to handle the term $\ln{\left[\frac{T}{m}(t^{1-\delta_3 }-1) \right]}$; in the way we find that the above term is of order
\begin{equation*}
O \left(  t^{-\frac{\delta_3 }{2}}  t^{\frac{1}{3} - \frac{4}{3}\sigma} \right), \quad \frac{1}{2} \le \sigma<1.
\end{equation*}

The error terms involving $1-t^{\delta_3  - 1}$ and $t^{-\delta_3 }$ are clearly smaller than the highest order term.

The second single sum of \eqref{sixthirtyeight-new} gives the following contribution:
\begin{equation*}
\frac{1}{t^{\frac{\delta_3  }{2}}} \frac{1}{T^{\sigma}} \sum_{m\in\tilde{N}} m^{-\sigma} m^{i (t-t^{\delta_3 })} \frac{1}{ \ln{ \left[  \frac{m}{T} (t^{1-\delta_3 }-1) \right]  }} \times O \left( 1, 1 - t^{\delta_3  - 1}, t^{-\delta_3 } \right) .
\end{equation*}
Since $m\in\tilde{N}$, we get the following bounds
\begin{equation} \label{8.40a}
\frac{1}{2\ln{t}} < \frac{1}{ \ln \left[  \frac{m}{T} \left( t^{1-\delta_3 }-1\right) \right] }  < t^\frac{\delta_3}{2},
\end{equation}
and the relevant sum can be analyzed in the same way as the single sum above. Thus, using \eqref{newsixfortyseven} it follows that the second sum above is of order $O(T^{1-2\sigma})$, or using the classical techniques it follows that this sum is of order $O( t^{ \frac{1}{3} -\frac{4}{3} \sigma } \ln{t})$, $\frac{1}{2}\le\sigma<1$.

\begin{figure}
\begin{center}
\includegraphics[scale=0.25]{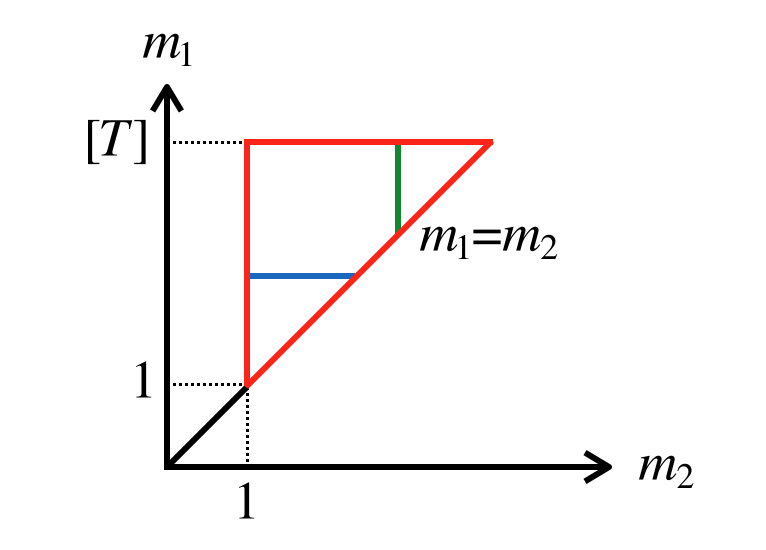}
\end{center}
\caption{Interchange of the order of summation.}
\label{Fig-s}
\end{figure}

Let $S_A$ denote the first error term involving a double sum:
$$S_A=\sum_{m_1=1}^{[T]}\sum_{m_2=1}^{[T]} \dfrac{1}{\left(m_1+m_2\right)^{\sigma-1+i t}}\dfrac{1}{m_2^{\sigma+1-i t}}.$$
Letting $m_2=m, \ m_1+m_2=n,$ and employing the triangular inequality we find
\begin{equation}\label{sixTsup1}
\big|S_A\big|=\Bigg|\sum_{m=1}^{[T]}\sum_{n=m+1}^{m+[T]} \dfrac{1}{n^{\sigma-1+i t}}\dfrac{1}{m^{\sigma+1-i t}}\Bigg| \leq \sum_{m=1}^{[T]} \Bigg| \sum_{n=m+1}^{m+[T]} \dfrac{1}{n^{\sigma-1+i t}}\Bigg| \dfrac{1}{m^{\sigma+1}}.
\end{equation}
Taking into consideration that $1\leq m \leq [T],$ and that $\sigma-1<0$ it follows that (see Appendix \ref{rem6}) 
\begin{equation}\label{sixTsup2}
\sum_{n=m+1}^{m+[T]} \dfrac{1}{n^{\sigma-1+i t}}=O\left(t^{\frac{3}{2}-\sigma}\right).
\end{equation}
Indeed, using \eqref{sixTsup2} into \eqref{sixTsup1} and noting that $\sigma+1>1$, it follows that 
\begin{equation}\label{sixTsup3}
\big|S_A\big|=O\left(t^{\frac{3}{2}-\sigma}\right).
\end{equation}
Let $S_B$ denote the second error involving a double sum:
\begin{equation}\label{sixTsup4}
S_B=\sum_{m_1=1}^{[T]}\sum_{m_2=1}^{[T]} \dfrac{1}{\left(m_1+m_2\right)^{\sigma-1+i t}}\dfrac{1}{m_2^{\sigma-i t}}\dfrac{1}{m_1}.
\end{equation}
Splitting this sum into two sums, depending on whether $m_1/m_2>1$ or $m_1/m_2<1,$ we find
\begin{equation}\label{sixTsup5}
S_B=S_1+S_2,
\end{equation}
where
\begin{equation}\label{sixTsup6}
S_1=\sum_{m_1=1}^{[T]}\sum_{m_2=1}^{m_1} \dfrac{1}{\left(m_1+m_2\right)^{\sigma-1+i t}}\dfrac{1}{m_2^{\sigma-i t}}\dfrac{1}{m_1},
\end{equation}
and
\begin{equation}\label{sixTsup7}
S_2=\sum_{m_1=1}^{[T]}\sum_{m_2=m_1+1}^{[T]} \dfrac{1}{\left(m_1+m_2\right)^{\sigma-1+i t}}\dfrac{1}{m_2^{\sigma-i t}}\dfrac{1}{m_1}.
\end{equation}
In order to estimate the sum $S_1$, we interchange the order of summation, see figure \ref{Fig-s}.

Thus,
\begin{equation*}
S_1=\sum_{m_2=1}^{[T]}\sum_{m_1=m_2}^{[T]} \dfrac{1}{\left(m_1+m_2\right)^{\sigma-1+i t}}\dfrac{1}{m_2^{\sigma-i t}}\dfrac{1}{m_1},
\end{equation*}
or
\begin{equation}\label{sixTsup8}
S_1=\sum_{m_2=1}^{[T]}\sum_{m_1=m_2}^{[T]} \dfrac{1}{\left(m_1+m_2\right)^{\sigma-1+i t}}\dfrac{1}{m_2^{\sigma+1-i t}}\dfrac{m_2}{m_1}.
\end{equation}
Using partial summation and the fact that $m_2/m_1\leq 1,$ it follows that
\begin{equation}\label{sixTsup9}
S_1=O\left(\tilde{S}_1\right),
\end{equation}
where
\begin{equation}\label{sixTsup10}
\tilde{S}_1=\sum_{m_2=1}^{[T]}\sum_{m_1=m_2}^{[T]} \dfrac{1}{\left(m_1+m_2\right)^{\sigma-1+i t}}\dfrac{1}{m_2^{\sigma+1-i t}}.
\end{equation}
Then, proceeding as with the sum $S_A$, it follows that
\begin{equation}\label{sixTsup11}
S_1=O\left(t^{\frac{3}{2}-\sigma}\right).
\end{equation}
In order to estimate $S_2$, we first note that
\begin{equation}\label{sixTsup12}
\big|S_2\big|\leq \sum_{m_1=1}^{[T]} \Bigg| \sum_{m_2=m_1+1}^{[T]} \dfrac{1}{\left(m_1+m_2\right)^{\sigma-1+i t}}\dfrac{1}{m_2^{\sigma-i t}}\Bigg| \dfrac{1}{m_1}.
\end{equation}
Then, taking into consideration that $m_1<m_2$, we can use the following ``crude" estimate for the $m_2$ sum:
\begin{equation}\label{sixTsup13}
\Bigg| \sum_{m_2=m_1+1}^{[T]} \dfrac{1}{\left(m_1+m_2\right)^{\sigma-1+i t}}\dfrac{1}{m_2^{\sigma-i t}}\Bigg| \leq \int_{m_1+1}^T \dfrac{1}{\left(m_1+x\right)^{\sigma-1}}\dfrac{1}{x^{\sigma}} dx:=J\left(m_1,t\right).
\end{equation}
But, $$m_1<x, \ \text{ or } \quad m_1+x<2x,  \ \text{ or } \quad \left(m_1+x\right)^{1-\sigma} <(2x)^{1-\sigma}.$$
Thus, $$J\left(m_1,t\right)< \int_{m_1+1}^T 2^{1-\sigma} x^{1-2\sigma} dx= O \left(T^{2-2\sigma}\right) + O \left(m_1^{2-2\sigma}\right)= O \left(t^{2-2\sigma}\right).$$
Hence, equation \eqref{sixTsup12} implies
\begin{equation}\label{sixTsup14}
\big|S_2\big| = O \left(t^{2-2\sigma}\int_{1}^T \frac{dx}{x}\right)= O \left(t^{2-2\sigma}\ln t\right).
\end{equation}

\textbf{QED}

\vphantom{a}

\begin{remark}\label{rem6.1}
Considering the case that the stationary points are on the end-point of the interval of integration, we get the following conditions $$\tau_1=\dfrac{1}{1+\lambda}=\dfrac{t^{\delta_2}}{t} \qquad or \qquad \tau_1=\dfrac{1}{1+\lambda}=1-\dfrac{t^{\delta_3}}{t},$$ equivalently
\begin{equation}\label{sixnr1}
\lambda=t^{1-\delta_2}-1\qquad or \qquad \lambda=\dfrac{1}{t^{1-\delta_3}-1}.
\end{equation}

In this case, the relevant contribution is computed using exactly the same procedure, but now we get half of the contribution, thus the rhs of \eqref{sixthirtyone} is multiplied by the term $\frac{1}{2}$.

The fact that the set of summation is restricted by the constraints \eqref{sixnr1}, makes the contribution of this sum negligible. In particular, this contribution is absorbed in the analysis of section \ref{sec7}. Indeed, if we have additional stationary points which are on the end-point of the interval of integration, then the strict inequalities that define the set $M$ in \eqref{sevenforteen} now allow equality. Thus, in this case the set $M$ has more elements, however, the remaining terms, namely the last two terms of the rhs of \eqref{sevenfifteen}, remain the same.

\end{remark}

Theorem \ref{t6.1} shows that the main contribution of the stationary points to $I_{3}$ is given by the double 
sum in (\ref{sixfortyseven}) involving $m_{2}^{-\bar{s} }(m_{1}+m_{2})^{-s}$. 
In what follows we present a rigorous analysis of this sum.

\section{The Analysis of the Double Sum Arising from the Stationary Points} \label{sec7}

In what follows we will compute the large $t$ asymptotics of the last two terms of the rhs of \eqref{sevensix}. 

\begin{lemma} \label{l7.2}
\begin{equation}
- \sum_{m=1}^{[T]}\frac{1}{m^{2\sigma }} + 2\Re\left\{ \sum_{m=1}^{[T]}\sum_{n=[T]+1}^{[T]+m}\frac{1}{m^{\bar{s}}n^{s} }  \right\}  = -
\frac{[T]^{1-2\sigma }}{1-2\sigma }+O(1)
\nonumber
\end{equation}
\begin{equation}
-\frac{1}{\pi }\Im \left\{ \frac{1}{T^{s}}
[\sum_{m=1}^{[T]}\frac{1}{m^{2\sigma -1}}
\bigl( \frac{1}{T}+\frac{1}{m} \bigr)^{1-s}\bigl( 1+O(\frac{1}{t}) \bigr)
+\bigl( \sum_{m=1}^{[T]}\frac{1}{m^{\bar{s} }}  \bigr) O(1)]   \right\},  \nonumber
\end{equation}
\begin{equation}  \label{sevenseven}
\hphantom{2a} s=\sigma +it, \hphantom{2a} 0<\sigma <1, \hphantom{2a} \sigma\neq\frac{1}{2}, \hphantom{2a} t\rightarrow \infty ,
\end{equation}
and
\begin{equation}
2\Re\left\{ \sum_{m=1}^{[T]}\sum_{n=[T]+1}^{[T]+m}\frac{1}{m^{\frac{1}{2}-it}n^{\frac{1}{2}+it}}  \right\} - \sum_{m=1}^{[T]}\frac{1}{m}
= - \ln t+O(1) \nonumber 
\end{equation}
\begin{equation}
-\frac{1}{\pi}\Im\left\{
\frac{1}{T^{\frac{1}{2}+it}}
[ \sum_{m=1}^{[T]}\bigl( \frac{1}{T}+\frac{1}{m} \bigr)^{\frac{1}{2}-it}\bigl( 1+O(\frac{1}{t})\bigr)
+\bigl( \sum_{m=1}^{[T]}\frac{1}{m^{\frac{1}{2}-it}}  \bigr)O(1) ]                
\right\},   \nonumber
\end{equation}
\begin{equation}  \label{seveneight}
t\rightarrow \infty .
\end{equation}
\end{lemma}

\textbf{Proof} 
We first show that
\begin{subequations} \label{sevennine}
\begin{equation} \label{sevenninea}
\sum_{m=1}^{[T]}\frac{1}{m^{2\sigma }}=\frac{[T]^{1-2\sigma }}{1-2\sigma }+O(1), \quad 0<\sigma <1, \quad \sigma \neq\frac{1}{2}, \quad t\rightarrow \infty,
\end{equation}
and
\begin{equation} \label{sevennineb}
\sum_{m=1}^{[T]}\frac{1}{m}=\ln t+O(1), \quad t\rightarrow \infty.
\end{equation}
\end{subequations}
To derive these formulae we will use the following identity (Theorem 2.1 of [T]):
\begin{equation}
\sum_{a < n\leq b}f(n)=\int_{a}^{b}f(x)dx+\int_{a}^{b}(x-[x]-\frac{1}{2})\frac{df(x)}{dx}dx \nonumber 
\end{equation}
\begin{equation}  \label{seventen}
+(a-[a]-\frac{1}{2})f(a)-(b-[b]-\frac{1}{2})f(b).  
\end{equation}
Letting in (\ref{seventen}) 
$$a=1,  \hphantom{2a}  b=[T], \hphantom{2a} f(x)=x^{-2\sigma },$$
and employing the identity
\begin{equation}   \label{seveneleven}
x-[x]-\frac{1}{2}=-\sum_{n=1}^{\infty }\frac{\sin (2n\pi x)}{n\pi }, 
\end{equation}
we find
\begin{equation}
\sum_{2}^{[T]}\frac{1}{m^{2\sigma }}=\int_{1}^{[T]}x^{-2\sigma }dx+O\bigl( \frac{1}{T^{2\sigma }} \bigr)+O(1)+
2\sigma \int_{1}^{[T]}\sum_{n=1}^{\infty}\frac{\sin (2n\pi x)}{n\pi}\frac{1}{x^{2\sigma +1}}.   \nonumber
\end{equation}
Thus,
\begin{equation}
\sum_{1}^{[T]}\frac{1}{m^{2\sigma }}=\int_{1}^{[T]}x^{-2\sigma }dx+O(1),  \nonumber
\end{equation}
and then equations (\ref{sevennine}) follow.

In order to simplify the double sum appearing in the lhs of equation (\ref{sevenseven}) we recall equation (1.10) of \cite{FL}:
\begin{equation}  \label{seventwelve}
\sum_{n=[T]+1}^{\left[\frac{\eta}{2\pi }\right]}\frac{1}{n^{s}}=\frac{1}{1-s}(\frac{\eta }{2\pi})^{1-s}+O(\frac{1}{t^{\sigma }}), 
\hphantom{2a} s=\sigma +it, \hphantom{2a} (1+\varepsilon ) t<\eta <\infty , \nonumber
\end{equation}
\begin{equation}
0\leq\sigma<1, \quad \varepsilon >0, \quad t\to \infty , 
\end{equation}
which is valid uniformly with respect to $\eta $ and $\sigma $.
Taking $\eta /2\pi =[t/2\pi ]+m$, equation (\ref{seventwelve}) becomes
\begin{equation}
\sum_{n=[T]+1}^{[T]+m}\frac{1}{n^{s}}=\frac{1}{1-s}(T+m)^{1-s}+O\bigl( \frac{1}{t^{\sigma }} \bigr)  \nonumber
\end{equation}
\begin{equation}
=\frac{i}{2\pi}\frac{1}{1+\frac{i(1-\sigma )}{t}}\frac{1}{T^{s}m^{s-1}}\bigl( \frac{1}{T}+\frac{1}{m} \bigr)^{1-s}
+O\bigl( \frac{1}{t^{\sigma }} \bigr) .  \nonumber
\end{equation}
Replacing in \eqref{sevenseven} the sum over $n$ by the above sum we find 
\begin{equation}
2\Re\left\{ \sum_{m=1}^{[T]}\sum_{n=[T]+1}^{[T]+m}\frac{1}{m^{\bar{s} }n^{s}} \right\}=
-\frac{1}{\pi }\Im\Biggl( 
\frac{1}{T^{s}}\sum_{m=1}^{[T]}\frac{1}{m^{2\sigma -1}}\bigl( \frac{1}{T}+\frac{1}{m}  \bigr)^{1-s}
\Bigl(1+O\bigl( \frac{1}{t} \bigr)\Bigr)   \nonumber
\end{equation}
\begin{equation}   \label{seventhirteen}
+O\bigl( \frac{1}{t^{\sigma }} \bigr)\sum_{m=1}^{[T]}\frac{1}{m^{\bar{s} }}\Biggr) , \hphantom{3a} t\rightarrow \infty .
\end{equation}
Equations (\ref{sevennine}) and (\ref{seventhirteen}) imply (\ref{seveneight}).

\textbf{QED}

\begin{figure}
\begin{center}
\includegraphics[width=\textwidth]{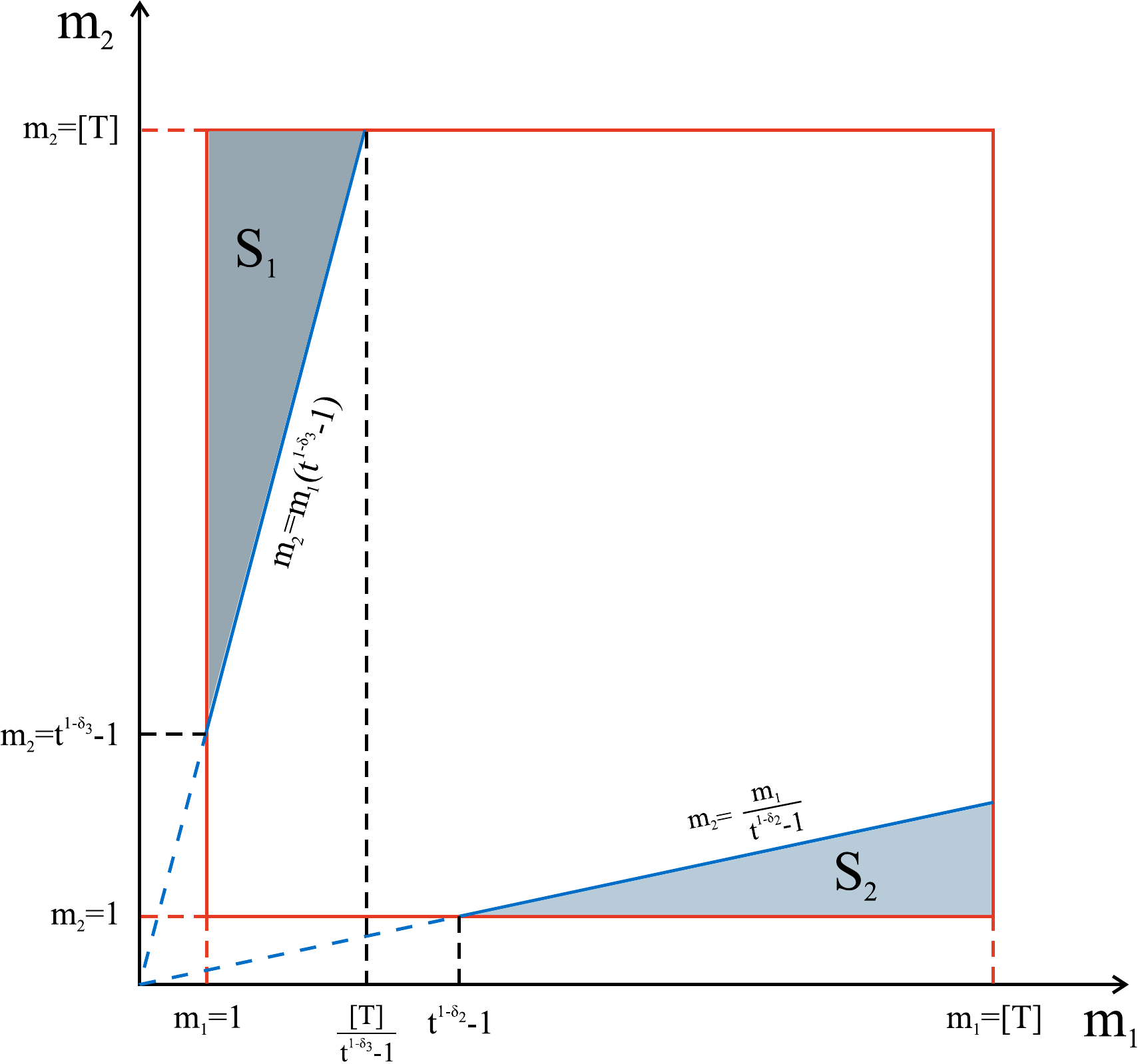}
\end{center}
\caption{The relation between the integers involved in the double sum without the constraint and the sum with the constraint.}
\label{fig2}
\end{figure}

\begin{remark}\label{rem7.1}
The single sums in \eqref{sevenseven} and \eqref{seveneight} involve the function $f(t,m)$ defined in \eqref{sixfortyfive}, thus the analysis of these sums is similar with the analysis of the Riemann zeta function via the techniques developed for the study of exponential sums. However, the occurrence of $\ t^{-\sigma}$ implies that it is sufficient to employ the following crude estimate:
\begin{equation}
\mid \sum_{m=1}^{[T]}\frac{1}{T^{s}m^{2\sigma -1}}\bigl( \frac{1}{T}+\frac{1}{m}  \bigr)^{1-s}\mid \leq
\int_{1}^{[T]}\frac{1}{T^{\sigma }x^{2\sigma -1}}\Bigl( \frac{1}{T}+\frac{1}{x}  \Bigr)^{1-\sigma }dx \nonumber
\end{equation}
\begin{equation}
=T^{1-2\sigma }\int_{1}^{T}\frac{1}{(\frac{x}{T})^{2\sigma -1}}\bigl( 1+\frac{T}{x}  \bigr)^{1-\sigma }
d\bigl( \frac{x}{T} \bigr) =
T^{1-2\sigma }\int_{\frac{1}{T}}^{1}\frac{( 1+\rho) ^{1-\sigma } }{\rho ^{\sigma}}d\rho   \nonumber
\end{equation}
\begin{equation}
= T^{1-2\sigma }O(1), \quad 0<\sigma<1.  \nonumber
\end{equation}
Indeed, $(1+\rho)^{1-\sigma}$ is bounded, thus,
\begin{equation*}
O \left( \int_{\frac{1}{T}}^{1}\frac{( 1+\rho) ^{1-\sigma } }{\rho ^{\sigma}}d\rho  \right) = O \left( \int_{\frac{1}{T}}^{1}\frac{d\rho}{\rho ^{\sigma}}  \right) = O(1).
\end{equation*}
\end{remark}

The above results show that the double sum of $m_{1}^{-s}(m_{1}+m_{2})^{-\bar{s} }$ appearing in the asymptotic evaluation of $I_{3}$ is related with the double sum expressing the large $t$ asymptotics of $|\zeta (s)|^{2}$. However, the former sum satisfies the constraint specified by the first inequality in (\ref{sixthirtynine}). In what follows we relate the sum without the constraint with the sum satisfying the above constraint, see figure \ref{fig2}.

\begin{lemma} \label{l7.3}
Define $M$ by 
\begin{multline}  \label{sevenforteen}
M = \left\{ m_{1}\in\mathbb{N}^{+}, \quad  m_{2}\in\mathbb{N}^{+}, \quad 1\leq m_{1}\leq [T], \quad 1\leq m_{2}< [T], \right.  \\
\left. \frac{1}{t^{1-\delta _{2}}-1}<\frac{m_{2}}{m_{1}}<t^{1-\delta _{3}}-1, \quad  t>0, \quad  T=\frac{t}{2\pi} \right\} ,
\end{multline}
where $\delta _{2}$ and $\delta _{3}$ are sufficiently small, positive constants.
Then 
\begin{equation}   \label{sevenfifteen}
\sum_{m_{1}=1}^{[T]}\sum_{m_{2}=1}^{[T]}=\sum_{m_{1},m_{2}\in M}+
\sum_{m_{1}=1}^{\frac{[T]}{t^{1-\delta _{3}}-1}-1}
\sum_{m_{2}=(t^{1-\delta _{3}}-1)m_{1}+1}^{[T]}+
\sum_{m_{1}=t^{1-\delta _{2}}}^{[T]}\sum_{m_{2}=1}^{\frac{m_{1}}{t^{1-\delta _{2}}-1}-1}.
\end{equation}
\end{lemma}

\textbf{Proof} 
\begin{equation}
\sum_{m_{1}=1}^{[T]}\sum_{m_{2}=1}^{[T]}=
\Biggl(  
\sum_{m_{1}=1}^{\frac{[T]}{t^{1-\delta _{3}}-1}-1}
+\sum_{m_{1}=\frac{[T]}{t^{1-\delta _{3}}-1}}^{t^{1-\delta _{2}}-1}
+\sum_{m_{1}=t^{1-\delta _{2}}}^{[T]}
\Biggr)\sum_{m_{2}=1}^{[T]}  \nonumber
\end{equation}
\begin{equation}   \label{sevensixteen}
=\sum_{m_{1}=1}^{\frac{[T]}{t^{1-\delta _{3}}-1}-1}\sum_{m_{2}=1}^{[T]}
+\sum_{m_{1}=\frac{[T]}{t^{1-\delta _{3}}-1}}^{t^{1-\delta _{2}}-1}\sum_{m_{2}=1}^{[T]}
+\sum_{m_{1}=t^{1-\delta _{2}}}^{[T]}\sum_{m_{2}=1}^{[T]}.
\end{equation}
We subdivide the sum over $m_{2}$ occuring in the first and third double sums in (\ref{sevensixteen}) as follows:

\begin{equation}
\sum_{m_{2}=1}^{[T]}=\sum_{m_{2}=1}^{(t^{1-\delta _{3}}-1)m_{1}}+
\sum_{m_{2}=(t^{1-\delta _{3}}-1)m_{1}+1}^{[T]}, \nonumber
\end{equation}
and
\begin{equation}
\sum_{m_{2}=1}^{[T]}= \sum_{m_{2}=1}^{\frac{m_{1}}{t^{1-\delta _{2}}-1}-1}
+\sum_{m_{2}=\frac{m_{1}}{t^{1-\delta _{2}}-1}}^{[T]}.  \nonumber
\end{equation}
Substituting the above expressions in (\ref{sevensixteen}) we find
\begin{equation}
\sum_{m_{1}=1}^{[T]}\sum_{m_{2}=1}^{[T]}=
\sum_{m_{1}=1}^{\frac{[T]}{t^{1-\delta _{3}}-1}-1}
\Biggl(
\sum_{m_{2}=1}^{(t^{1-\delta _{3}}-1)m_{1}}+
\sum_{m_{2}=(t^{1-\delta _{3}}-1)m_{1}+1}^{[T]}
\Biggr)   \nonumber
\end{equation}
\begin{equation}  \label{sevenseventeen}
+\sum_{m_{1}=\frac{[T]}{t^{1-\delta _{3}}-1}}^{t^{1-\delta _{2}}-1}\sum_{m_{2}=1}^{[T]}+
\sum_{m_{1}=t^{1-\delta _{2}}}^{[T]}
\Biggl(
\sum_{m_{2}=1}^{\frac{m_{1}}{t^{1-\delta _{2}}-1}-1}+
\sum_{m_{2}=\frac{m_{1}}{t^{1-\delta _{2}}-1}}^{[T]}
\Biggr).
\end{equation}
The sum of the first, third, and fifth double sums in (\ref{sevenseventeen}) equals the first term of the rhs of 
(\ref{sevenfifteen}), whereas the second and fourth double sums in (\ref{sevenseventeen}) are the second and third 
terms of the rhs of (\ref{sevenfifteen}).

\textbf{QED}

\vphantom{a}


Equation \eqref{sevenfifteen} shows that the second double sum occurring in the asymptotic evaluation of $I_3$ in \eqref{sixfortyseven} differs from the double sum occuring in the asymptotics of $|\zeta(s)|^2$ by the two sums defined in the second and third terms of the rhs of \eqref{sevenfifteen}. By combining certain results of \cite{FL} with techniques developed in \cite{T2} it is possible to estimate the latter two sums.

\begin{lemma} \label{lemma7.4}
Define the double sum $S_{1}$ by 
\begin{equation*}
S_{1}(\sigma,t,\delta )=\sum_{m_{1}=1}^{\frac{[T]}{t^{1-\delta }-1}-1}\sum_{m_{2}=(t^{1-\delta }-1)m_{1}+1}^{[T]}
\frac{1}{m_{1}^{s}(m_{1}+m_{2})^{\bar{s} }},
\end{equation*}
\begin{equation} \label{seveneighteen}
s=\sigma +it, \quad 0<\sigma <1, \quad t>0,
\end{equation}
and $\delta $ is a sufficiently small, positive constant. Then,
\begin{equation} \label{sevennineteen}
S_{1}(\sigma,t,\delta ) = O\left( t^{\frac{1}{2}-\sigma}\tilde{G}(\sigma,t,\delta) \right) + O\left( \frac{t^{(1-\sigma)\delta}}{t^\sigma} \right), ~ 0<\sigma<1, ~ t\to\infty,
\end{equation}
where
\begin{equation} \label{seventwenty}
\tilde{G}(\sigma,t,\delta ) = O \left(  t^{(1-\sigma)\delta} \right) + O \left(  t^{\sigma\delta} \right), \quad 0<\sigma<1, \quad \sigma\ne\frac{1}{2}, \quad t\to\infty,
\end{equation}
and
\begin{equation} \label{seventwentyone}
\tilde{G}\left(\frac{1}{2},t,\delta\right) = O \left(  t^{\frac{\delta}{2}} \ln{t} \right), \quad t\to\infty.
\end{equation}
\end{lemma}
\textbf{Proof}
Letting $m_1=m$ and $m_1+m_2=n$ in the definition \eqref{seveneighteen} of $S_1$ yields
\begin{equation} \label{seventwentytwo}
S_{1}(\sigma,t,\delta )=\sum_{m=1}^{\frac{T}{t^{1-\delta }-1}-1}\sum_{n=t^{1-\delta }m+1}^{[T]+m}
\frac{1}{m^{s}n^{\bar{s}}}.
\end{equation}
The upper limit of the $m$-sum is given by
\begin{equation} \label{seventwentythree}
\frac{T}{t^{1-\delta }-1}-1 = \frac{t}{2\pi t^{1-\delta }\left(1-t^{\delta -1}\right)}-1 = \frac{t^\delta}{2\pi}L(t),
\end{equation}
where $L(t)$ is defined by
\begin{equation} \label{seventwentyfour}
L(t) = \frac{1}{1-t^{\delta - 1}} - \frac{2\pi}{t^\delta}.
\end{equation}
Hence,
\begin{equation} \label{seventwentyfive}
L(t) =1 - \frac{2\pi}{t^\delta} + O \left( \frac{t^\delta}{t}  \right), \quad t\to\infty.
\end{equation}
If $m=1$ then $t^{1-\delta}m=t^{1-\delta}$, and if $m=t^{\delta}L(t)/2\pi$, then $t^{1-\delta}m=tL(t)/2\pi$. Thus,
\begin{equation} \label{seventwentysix}
t^{1-\delta} \leq n \leq \frac{t}{2\pi} L(t).
\end{equation}
It is convenient to split the $S_1$ sum in terms of the following two sums:
\begin{equation} \label{seventwentyseven}
S_{A}(\sigma,t,\delta )=\sum_{m=1}^{\frac{t^{\delta}}{2\pi}L(t)}\sum_{n=t^{1-\delta }m+1}^{\left[\frac{t}{2\pi}\right]}
\frac{1}{m^{s}n^{\bar{s}}}, \quad 0<\sigma <1, \quad t>0,
\end{equation}
and
\begin{equation} \label{seventwentyeight}
S_{B}(\sigma,t,\delta )=\sum_{m=1}^{\frac{t^{\delta}}{2\pi}L(t)}\sum_{n=\left[\frac{t}{2\pi}\right]+1}^{\left[\frac{t}{2\pi}\right]+m}
\frac{1}{m^{s}n^{\bar{s}}}, \quad 0<\sigma <1, \quad t>0.
\end{equation}
Thus, computing $S_1$ reduces to computing $S_A$ and $S_B$:
\begin{equation} \label{seventwentynine}
S_{1}(\sigma,t,\delta ) = S_{A}(\sigma,t,\delta ) + S_{B}(\sigma,t,\delta ).
\end{equation}

For the rigorous  derivation of the estimates of $S_A$ and $S_B$, we refer to the proof of Theorem 5.1 in \cite{FK}, where the following results are obtained:

\begin{equation} \label{sevenfortytwo}
S_A = O \left( t^{\frac{1}{2}-\sigma}\right) \tilde{G}(\sigma,t, \delta), \quad t\to\infty,
\end{equation}
and
\begin{equation} \label{seventhirtyfour}
S_{B}(\sigma,t,\delta) = O\left( t^{-\sigma+(1-\sigma)\delta} \right), \quad t\to\infty,
\end{equation}
where $\tilde{G}(\sigma,t, \delta)$ is given by \eqref{seventwenty} and \eqref{seventwentyone}.

Equations \eqref{seventwentynine}, \eqref{sevenfortytwo} and \eqref{seventhirtyfour} imply \eqref{sevennineteen}.

\textbf{QED}

\begin{lemma} \label{l7.5}
Define the double sum $S_2(\sigma,t,\delta)$ by
\begin{equation} \label{sevenfortythree}
S_{2}(\sigma,t,\delta )=\sum_{m_{1}=t^{1-\delta}}^{[T]} \sum_{m_{2}=1}^{\frac{m_1}{t^{1-\delta}-1}-1}
\frac{1}{m_{1}^{s}(m_{1}+m_{2})^{\bar{s} }}, ~ s=\sigma+it, ~ 0<\sigma<1, ~t>0,
\end{equation}
where $\delta$ is a sufficiently small, positive constant. Then,
\begin{equation} \label{sevenfortyfour}
S_{2}(\sigma,t,\delta ) = O \left(  t^{1-2\sigma+2\delta\sigma} (\ln{t})^3 \right) + O \left(  \frac{t^\delta}{t^{2\sigma} }\right), \quad 0<\sigma<1, \quad t\to\infty.
\end{equation}
\end{lemma}
\textbf{Proof}
Letting $m_1=m$ and $m_1+m_2=n$ in the definition \eqref{sevenfortythree} of $S_2$ we find
\begin{equation} \label{sevenfortyfive}
S_{2}(\sigma,t,\delta ) = \sum_{m=t^{1-\delta}}^{[T]} \sum_{n=1+m}^{\frac{m}{1-t^{\delta-1}}-1}
\frac{1}{m^{s}n^{\bar{s} }}.
\end{equation}
It is convenient to split the $S_2$ sum in terms of the following two sums:
\begin{equation} \label{sevenfortysix}
S_{A}(\sigma,t,\delta )=\sum_{m=t^{1-\delta}}^{[T]}\sum_{n=1+m}^{[T]}
\frac{1}{m^{s}n^{\bar{s}}}, \quad 0<\sigma <1, \quad t>0,
\end{equation}
and
\begin{equation} \label{sevenfortyseven}
S_{B}(\sigma,t,\delta )=\sum_{m=t^{1-\delta}}^{[T]}\sum_{n=[T]+1}^{\frac{m}{1-t^{\delta-1}}-1}
\frac{1}{m^{s}n^{\bar{s}}}, \quad 0<\sigma <1, \quad t>0.
\end{equation}
Hence
\begin{equation} \label{sevenfortyeight}
S_{2}(\sigma,t,\delta) = S_{A}(\sigma,t,\delta) + S_{B}(\sigma,t,\delta), \quad 0<\sigma <1, \quad t>0.
\end{equation}
We first analyze $S_B$. In the proof of Theorem 5.2 in \cite{FK} the following estimate is obtained:
\begin{equation} \label{sevenfiftythree}
S_B(\sigma,t,\delta) = O\left(  \frac{t^{\delta}}{t^{2\sigma}} \right), \quad 0<\sigma<1, \quad t \to\infty.
\end{equation}

We next consider $S_A$. This sum can be analyzed in two different ways:

The first involves changing the order of summation, employing the asymptotic formula (2.7) in \cite{FK}, and proceeding as in the proof of $S_A$ in Lemma \ref{lemma7.4}. A rigorous proof for this estimate is given in \cite{FK}.

 The second way uses the techniques developed in \cite{T2} and \cite{K} and appropriately modified in Appendix \ref{appDS}. These techniques imply that
\begin{equation} \label{sevenfiftysix}
\sum_{m=t^{1-\delta}}^{[T]} \sum_{n=m+1}^{[T]} \frac{1}{m^{it}n^{-it}} = O \left(  t(\ln{t})^3 \right), \quad t\to \infty.
\end{equation}
Furthermore, under the condition that the expressions 
\begin{equation}\label{psks}
b_{m,n} - b_{m+1,n}, \quad b_{m,n}-b_{m,n+1}, \quad b_{m,n}-b_{m+1,n}-b_{m,n+1}+b_{m+1,n+1},\end{equation}  
keep their sign, the following result is derived in \cite{T2}:
\begin{equation} \label{sevenfiftyfive1}
\left|  \sum_{m=1}^{M} \sum_{n=1}^{N} a_{m,n} b_{m,n} \right| \le 5GH,
\end{equation}
where
\begin{equation} \label{sevenfiftysix1}
S_{m,n} \Doteq \sum_{\mu=1}^{m} \sum_{\nu=1}^{n} a_{\mu,\nu}, \quad |S_{m,n}|\le G, \quad 1\le m\le M, \quad 1\le n\le N,
\end{equation}
with
\begin{equation} \label{sevenfiftyseven}
b_{m,n} \in \mathbb{R}, \quad 0 \le b_{m,n} \le H.
\end{equation}
We apply the above argument for $\displaystyle b_{m,n}=\frac{1}{m^\sigma n^\sigma}$,
thus  the expressions in \eqref{psks} keep their sign, and furthermore $$H=\dfrac{1}{t^{(1-\delta)\sigma}}\dfrac{1}{t^{(1-\delta)\sigma}}=t^{-2\sigma}t^{2\delta\sigma}.$$

Combining the above result with \eqref{sevenfiftysix} we find
\begin{equation} \label{sevenfiftyeight}
S_A(\sigma,t,\delta) = O \left( t^{1-2\sigma} t^{2\delta\sigma} (\ln{t})^3 \right), \quad 0<\sigma<1, \quad t \to\infty.
\end{equation}


Equations \eqref{sevenfortyeight}, \eqref{sevenfiftythree}, \eqref{sevenfiftyeight} imply \eqref{sevenfortyfour}.

\textbf{QED}

\vphantom{a}

Combining equation \eqref{sevensix} with Lemmas \ref{l7.2}, \ref{l7.3}, \ref{lemma7.4}, \ref{l7.5}, we obtain the following result:

\begin{theorem}\label{t7.1}
Let the set $M$ be defined in \eqref{sevenforteen}. Then,
\begin{multline} \label{sevenfiftynine}
2\Re \left\{ \sum_{m_1,m_2 \in M}\frac{1}{m_1^{s}\left(  m_1 + m_2 \right)^{\bar{s}}} \right\} - \left(  \sum_{m=1}^{[T]} \frac{1}{m^s} \right) \left(  \sum_{m=1}^{[T]} \frac{1}{m^{\bar{s}}} \right) = \\
A(\sigma,t) + O\left( t^{-\sigma} \right) \Im{\left\{  \sum_{m=1}^{[T]} \frac{1}{m^{\bar{s}}} \right\}} + O\left( t^{1-2\sigma+2\delta_2\sigma} (\ln{t})^3 \right) + O\left(  \frac{t^{\delta_2}}{t^{2\sigma}}  \right) + O(1) \\
- \begin{cases}
      \frac{[T]^{1-2\sigma}}{1-2\sigma} + O\left( t^{\frac{1}{2}-\sigma} \right) \left(  O\left(  t^{(1-\sigma)\delta_3}  \right) + O\left(  t^{\sigma \delta_3} \right)  \right) + O \left(  \frac{t^{(1-\sigma)\delta_3}}{t^{\sigma}} \right) , &  \sigma \ne \frac{1}{2}, \\
      \ln{t} + O\left(  t^{\frac{\delta_3}{2}} \ln{t} \right), & \sigma = \frac{1}{2},
      \end{cases}
\end{multline}
where
\begin{equation} \label{sevensixty}
A(\sigma, t) = - \frac{1}{\pi} \Im{ \left\{ \frac{1}{T^s} \sum_{m=1}^{[T]} \frac{1}{m^{2\sigma-1}} \left(  \frac{1}{T} +\frac{1}{m} \right)^{1-s} \left( 1 + O\left(  \frac{1}{t} \right)  \right)  \right\}}
\end{equation}
and
\begin{equation*}
A(\sigma, t) = O \left(  t^{1-2\sigma} \right), \quad t\to\infty.
\end{equation*}
\end{theorem}

\textbf{Proof} Equation \eqref{sevenfifteen} yields
\begin{multline} \label{sevensixtyone}
2\Re \left\{ \sum_{m_1,m_2 \in M} \frac{1}{m_1^s(m_1+m_2)^{\bar{s} }} \right\} = 2\Re \left\{ \sum_{m_1=1}^{[T]}\sum_{m_2=1}^{[T]}  \frac{1}{m_1^s(m_1+m_2)^{\bar{s} }} \right\} \\
-2 \Re{\left\{  S_1(\sigma, t, \delta_2) + S_2(\sigma, t, \delta_3)  \right\}}.
\end{multline}
Replacing the first term of the rhs of \eqref{sevensixty} with the rhs of \eqref{sevensix}, equation \eqref{sevensixtyone} becomes
\begin{equation} \label{sevensixtytwo}
\text{LHS} = - \sum_{m=1}^{[T]} \frac{1}{m^{2\sigma}} + 2\Re \left\{ \sum_{m=1}^{[T]}\sum_{n=[T]+1}^{[T]+m}  \frac{1}{m^{\bar{s} }n^s} \right\} - 2\Re \left\{ S_1(\sigma, t, \delta_3) + S_2(\sigma, t, \delta_2) \right\},
\end{equation}
where LHS denotes the lhs of \eqref{sevenfiftynine}.

According to equations \eqref{sevenseven} and \eqref{seveneight}, the first term of the rhs  of \eqref{sevensixtytwo} is given by
\begin{multline*}
A(\sigma,t) + O(t^{-\sigma}) \Im{\left\{ \sum_{m=1}^{[T]} \frac{1}{m^{\bar{s}}}  \right\}} + O(1) - \begin{cases}
      \frac{[T]^{1-2\sigma}}{1-2\sigma}, &  \sigma \ne \frac{1}{2}, \\
      \ln{t}, & \sigma = \frac{1}{2},
      \end{cases}, \quad t\to \infty,
\end{multline*}
where $A(\sigma,t)$ is defined in \eqref{sevensixty}. Replacing the first term of the rhs of \eqref{sevensixtytwo} by the above expression, as well as replacing $S_1$ and $S_2$ by the rhs of equations \eqref{sevennineteen} and \eqref{sevenfortyfour}, with $\delta$ replaced by $\delta_2$ and $\delta_3$, respectively, equation \eqref{sevensixtytwo} becomes equation \eqref{sevenfiftynine}.

\textbf{QED}

\begin{remark}\label{rem7.2}
The error term appearing in Theorem \ref{t7.1} can be simplified as follows:
\begin{multline} \label{sevensixtythree}
2\Re \left\{ \sum_{m_1,m_2 \in M}\frac{1}{m_1^{s}\left(  m_1 + m_2 \right)^{\bar{s}}} \right\} - \left(  \sum_{m=1}^{[T]} \frac{1}{m^s} \right) \left(  \sum_{m=1}^{[T]} \frac{1}{m^{\bar{s}}} \right) = \\
\begin{cases} O\left( t^{1+2\left(\delta_2-1\right) \sigma} \left( \ln t \right)^3\right), & 0<\sigma<\frac{1}{2}+O\left(\delta_2\right), \\
 O\left( t^{\delta_2} \left( \ln t \right)^3\right)  + O\left( t^{\frac{\delta_3}{2}}  \ln t \right), & \sigma=\frac{1}{2}, \\
O(1), & \frac{1}{2}+O\left(\delta_2\right) <\sigma<1.
\end{cases}
\end{multline}
Moreover, the first term of the second case is removed by choosing $\delta_3> 2 \delta_2$.
\end{remark}

\section{The main result and Conclusions}

\underline{\textbf{Proof of \eqref{1.14}}}

Equation \eqref{1.14} is a direct consequence of equations (1.9) and theorems \ref{t5.1}, \ref{t6.1} and \ref{t7.1}.
Regarding the expressions for $I_3$ and $I_4$ we note the following:
\begin{itemize}
\item[(i)] In Theorem \ref{t5.1} we choose $\delta=\delta_3> \delta_4,$ and also employ Remark \ref{rem5}.
\item[(ii)] In Theorem \ref{t6.1}, $\delta=\delta_3$.
\item[(iii)] In Theorem \ref{t7.1} we choose $\delta=\delta_3 > 2\delta_2,$ and also employ the formulae of Remark \ref{rem7.2}.
\end{itemize}
Then, substituting these expressions of $I_3$ and $I_4$ in equations \eqref{4.21} and \eqref{4.22}, which are given by theorems \ref{t5.1} and \ref{t6.1}, we obtain
\begin{multline}\label{8.41}
t^{-\frac{\delta}{2}} \frac{\sqrt{2}}{\pi}  \left[ 1 + O\left(   t^{-\delta}\right) \right]  \times \Re \Bigg\{ e^{\frac{i\pi}{4}} \left(  t^{\delta-1}\right)^{it^{\delta}} \left(  1 -  t^{\delta-1}\right)^{\sigma- \frac{1}{2}+i(t-t^{\delta})} \times \\
 \mathop{\sum\sum}_{m_1,m_2\in N} \frac{1}{m_1^{s-it^{\delta}}}  \frac{1}{m_2^{\bar{s} + it^{\delta}}} \frac{1}{\ln{\left[ \frac{m_2}{m_1} \left(  t^{1-\delta} - 1 \right)  \right]}} \big[ 1 + o(1)\big]  \Bigg\} \\
 = \begin{cases} -\ln t + O\left(t^\frac{\delta}{2} \ln t\right), & \sigma=\frac{1}{2}, \\
-\zeta(2\sigma) +O(1), & \frac{1}{2}+O(\delta)<\sigma<1, \end{cases}, \qquad t \to \infty. 
\end{multline}

Equation \eqref{1.14} follows from \eqref{8.41}.

\textbf{QED}

The main result of this paper is the rigorous derivation
of equation \eqref{1.14}, which provides the proof
of the analogue of Lindel\"{o}f's hypothesis for a slight variant
of $|\zeta(s)|^{2}$. Indeed,
the main difference of the double sum in the lhs of \eqref{1.14}  from the double sum characterising the leading 
large $t$-asymptotics of $|\zeta(s)|^{2}$ is the term $\ln [\frac{m_2}{m_1}(t^{1-\delta_{3}}-1)]$, and this term is larger than $\frac{1}{2\ln t}$ and  smaller than $t^{\frac{\delta_{3}}{2}}$.
Thus,
we expect that the sum in \eqref{1.14} times the term $\ln t$ behaves for large $t$ like $|\zeta(s)|^{2}$;
the possibility of establishing rigorously this result is discussed in \cite{FK}.

The derivation of \eqref{1.14} is split in two main parts. The first part is rather straightforward and it leads to the derivation of equations (1.9). Regarding this part we note that
the derivation of equation \eqref{1.3} is based on a certain identity
relating the Riemann and the Hurwitz zeta functions derived in \cite{AsF}, and on
the use of the Plemelj  formulae.
The proof that $I_{1}$ and $I_{2}$ are ``small" is based on  Atkinson's
classical estimates and on the use of the second mean value theorem for integrals.
The second part of the derivation of \eqref{1.14}  involves the derivation of the  rigorous estimation of the error terms in the asymptotic evaluation
of the integrals $I_{3}$ and $I_{4}$ . The main difficulty of this part is
the estimation of the terms arising from keeping
 the next order
term of the asymptotic expansion of $\zeta(s)$ (this term  is given in Theorem 3.1 of \cite{FL}). In this connection, a crucial role is played by the explicit asymptotic identities (5.1) and (5.3) of section 5
of \cite{FL}. Using these identities, together with some ``crude" estimates where
certain sums are replaced by integrals, it is   possible in most cases to obtain results which are at least as good as estimates obtained via the classical techniques for single and double sums. The reason for this  effectiveness of the identities (5.1) and (5.3) of  \cite{FL} is explained in \cite{FK}.

\begin{appendices}

\section{(proof of \eqref{sixTsup2}) }\label{rem6}
Let $\chi(s)$  be defined by
\begin{equation}\label{sixrem1}
\chi(s) = \frac{(2\pi)^s}{\pi} \sin{\left( \frac{\pi s}{2} \right)} \Gamma(1-s), \qquad s\in\mathbb{C}.
\end{equation}
It is shown in \cite{FL} that
\begin{equation}\label{sixrem2}
\chi(s)=\left(\frac{2\pi}{t}\right)^{s-\frac{1}{2}} e^{it} e^{i\frac{\pi}{4}}\left[1+O\left(\frac{1}{t}\right)\right], \qquad s=\sigma+it, \quad \sigma\in \mathbb{R},  \quad t\to \infty.
\end{equation}
Employing the well known identity
\begin{equation}\label{sixrem3}
\zeta(s)=\chi(s) \zeta(1-s), \qquad s\in \mathbb{C},
\end{equation}
with $s=\sigma-1+it,$ we find 
\begin{equation}\label{sixrem4}
\zeta(\sigma-1+it)=\chi(\sigma-1+it) \zeta(2-\sigma-it).
\end{equation}
Suppose that $0<\sigma<1$. Using the fact that $ \zeta(2-\sigma-it)$ is bounded as $t\to\infty$, as well as the asymptotic estimate \eqref{sixrem2}, equation \eqref{sixrem4} implies that
\begin{equation}\label{sixrem5}
\zeta(\sigma-1+it)=O\left(t^{\frac{3}{2}-\sigma}\right), \qquad 0<\sigma<1,\quad t\to\infty.
\end{equation}

The following result is derived in \cite{FL}, see Theorem 3.2, equation (3.20):
\begin{align}\label{sixrem6}
\zeta(s) = & \sum_{n =1}^{[\frac{t}{2\pi}]} n^{-s} - \frac{1}{1-s} \left(\frac{t}{2\pi}\right)^{1-s} 
	\\ \nonumber
& + \frac{e^{-\frac{i\pi (1-s)}{2}}}{(2\pi)^{1-s}}\sum_{n=1}^\infty  \sum_{j=0}^{N-1} e^{-nz - it\ln{z}} \left(\frac{1}{n + \frac{it}{z}}\frac{d}{dz}\right)^j\frac{z^{-\sigma}}{n + \frac{it}{z}} \biggr|_{z = it}
	\\ \nonumber
& + \frac{e^{\frac{i\pi (1-s)}{2}}}{(2\pi)^{1-s}}\sum_{n=2}^\infty
\sum_{j=0}^{N-1} e^{-nz - it\ln{z}} \left(\frac{1}{n + \frac{it}{z}}\frac{d}{dz}\right)^j\frac{z^{-\sigma}}{n + \frac{it}{z}} \biggr|_{z = -it}
	\\ \nonumber
& + \left(\frac{t}{2\pi}\right)^{1-s} e^{it}
\sum_{k=0}^{2N}   \frac{\overline{c_k(1-\sigma)}\Gamma(\frac{k+1}{2})}{t^{\frac{k+1}{2}}}
 +  O\biggl((2N +1)!! N 2^{2N} t^{-\sigma - N}  \biggr),
	\\ \nonumber
& \hspace{2cm}  0 \leq \sigma \leq 1, \quad N \geq 2, \quad t \to \infty,	
\end{align}
where the error term is uniform for all $\sigma, N$ in the above ranges and the coefficients $c_k(\sigma)$ are given therein.
This equation is derived in \cite{FL} under the assumption that $0<\sigma<1$. However, it is straightforward to verify that it is also valid for $-1<\sigma<0.$
Equations \eqref{sixrem6} and \eqref{sixrem5} imply that
\begin{equation}\label{sixrem7}
\sum_{m=1}^{[T]} \dfrac{1}{m^{\sigma-1+it}} = O\left(t^{\frac{3}{2}-\sigma}\right), \qquad 0<\sigma<1,\quad t\to\infty.
\end{equation}

%
%

\section{(proof of \eqref{sevenfiftysix})}\label{appDS}

We prove that
\begin{equation}\label{DS1}
\sum_{m=M}^{M'}\sum_{n=N}^{N'} \dfrac{1}{m^{it}} \dfrac{1}{n^{-it}} = O(t \ln t),
\end{equation}
with $n>m$, and $$\begin{cases}A_1\sqrt{t}<M<M'<2M<A_3 t, \\ A_2\sqrt{t}<N<N'<2N<A_4 t,\end{cases}$$ for some positive constants $\{A_j\}_1^4$.

In this connection, we divide the set of summation similarly to the division implemented in Theorem 1 of \cite{T2}, namely, in ``small" rectangles $\Delta_{p,q}$, such that 
 $$\begin{cases}  M+p \l_1\leq m \leq M+p\l_1+\l_1, \\ N+q \l_2\leq n \leq N+q\l_2+\l_2. \end{cases}$$ 
Moreover, we pick 
\begin{equation}\label{DS2}
\begin{cases} \l_1=c_1\frac{M^2}{t}, \\ \l_2=c_2\frac{N^2}{t}, \end{cases}
\end{equation} for some positive constants $c_1$ and $c_2$.

We make the following observations:
\begin{itemize}
\item $\displaystyle \begin{cases} 1 \leq \l_1 \leq M \Leftrightarrow  A_1\sqrt{t}<M<A_3 t, \\ 1 \leq \l_2 \leq N \Leftrightarrow  A_2\sqrt{t}<N<A_4 t,\end{cases}$ for some positive constants $\{A_j\}_1^4$.
\item The number of the ``small" rectangles $\Delta_{p,q}$  is $O\left(\dfrac{MN}{\l_1 \l_2}\right)$.
\end{itemize}
We use Theorem 2.16 of \cite{K} with $$f(x,y)=t(\ln x - \ln y).$$
Then, in each rectangle $\Delta_{p,q}$, with $n>m$ (equivalently $x>y$), the conditions of this theorem are satisfied with $\lambda_1=\frac{t}{M^2}$ and $\lambda_2=\frac{t}{N^2}$, because $$\big| f_{xx} \big| =\frac{t}{x^2}, \quad \big| f_{yy} \big| =\frac{t}{y^2} \quad \text{ and } \ \big| f_{xy} \big| =0.$$
Using the following facts: 
\begin{itemize}
\item the conditions  $\displaystyle \begin{cases} M > A_1 \sqrt{t}, \\ N> A_2 \sqrt{t}, \end{cases}$ imply that $\displaystyle \begin{cases} \lambda_1 <\frac{1}{A_1^2}, \\ \lambda_2 <\frac{1}{A_2^2} \end{cases}$, 
\item all the quantities $ \ln \big|\Delta_{p,q} \big|, \ \big|\ln \lambda_1 \big|$ and $\big|\ln \lambda_2 \big|$ are of order $O(\ln t)$,
\end{itemize}
and employing equation (2.56) of \cite{K}, we find
\begin{equation}\label{DS3}
\mathop{\sum\sum}_{(m,n)\in \Delta_{p,q}} e^{if(m,n)} = O\left(\dfrac{\ln t}{\frac{t}{MN}}\right).
\end{equation}
Thus, the fact that the number of the rectangles $\Delta_{p,q}$  is $O\left(\dfrac{MN}{\l_1 \l_2}\right)$, implies that
\begin{equation}\label{DS4}
\sum_{m=M}^{M'}\sum_{n=N}^{N'} e^{if(m,n)} = O\left(\dfrac{MN}{t} \ln t \dfrac{MN}{\l_1 \l_2}\right).
\end{equation}
Equation \eqref{DS1} follows from applying \eqref{DS2} in \eqref{DS4}.

Finally, using the classical splitting for the sets of summation for exponential sums, see \cite{T} and \cite{T2}, equation \eqref{sevenfiftysix} follows from applying \eqref{DS1} for  $O\left(\left(\ln t^\delta\right)^2\right)=O\left(\left(\delta\ln t\right)^2\right)=O\left((\ln t)^2\right)$ times.

\end{appendices}

\section*{Acknowledgment}

This project would {\it not} have been completed without the crucial contribution of Kostis Kalimeris. Kostis has studied extensively the classical techniques for the estimation of single and multiple exponential sums; these techniques are used extensively in our joint paper with Kostis \cite{FK} and some of the results of this paper are used in sections \ref{sec6} and \ref{sec7}. Furthermore, Kostis has checked the entire manuscript and has made important contributions to the completion of Theorems \ref{t6.1} and \ref{t7.1}.

I have benefited greatly from my long collaboration with Jonatan Lenells. In particular, regarding the current work, equation \eqref{1.3} was derived in June 2015 as part of a long term collaborative project with Jonatan on the asymptotics of the Riemann zeta function and of related functions; furthermore, several results from our joint paper with Jonatan \cite{FL} are used in the present paper.

Arran Fernandez and Euan Spence made important contributions regarding
the rigorous estimates of the term $I_{B}$.

The starting point of the approach developed here is equation (2.5) of \cite{F} which is derived in our joint paper with Anthony Ashton \cite{AsF}.

In addition to my former students Anthony, Euan and Kostis, my current student Arran, and my former post doctoral associate Jonatan, my former students Mihalis Dimakos and Dionysis Mantzavinos have offered me generous support and assistance during the last seven years of my investigation of the asymptotics of the Riemann zeta function.

I am grateful to the late Bryce McLeod, as well as to Eugene Shargorodsky and Bengt Fornberg for collaborative attempts related to the present paper.

I thank my current student Nicholas Protonotarios, as well as my current post doctoral associates Parham Hashemzadeh and Iason Hitzazis for their help.

Finally, I am deeply grateful to EPSRC for many years of continuous support which currently is in the form of a senior fellowship.

\end{document}